\documentclass[11pt, reqno]{amsart}

\usepackage{amsmath}
\usepackage{fullpage}
\usepackage{amssymb}
\usepackage{amsthm}
\usepackage{bbm}
\usepackage{subfigure}
\usepackage{float}
\usepackage{graphicx}
\usepackage{epstopdf}
\usepackage{color}
\graphicspath{{../figures/}}

\newtheorem{lemma}{Lemma}[section]
\newtheorem{theorem}[lemma]{Theorem}
\newtheorem{proposition}[lemma]{Proposition}
\newtheorem{remark}[lemma]{Remark}
\newtheorem{corollary}[lemma]{Corollary}
\title{Low regularity error estimates\\ for the time integration of 2D NLS}

\author{Lun Ji}
\address{Department of Mathematics, Universit\"{a}t Innsbruck, Technikerstr.~13, 6020 Innsbruck, Austria (L.~Ji)}
\email{lun.ji@uibk.ac.at}

\author{Alexander Ostermann}
\address{Department of Mathematics, Universit\"{a}t Innsbruck, Technikerstr.~13, 6020 Innsbruck, Austria (A.~Ostermann)}
\email{alexander.ostermann@uibk.ac.at}

\author{Fr\'ed\'eric Rousset}
\address{Universit\'e Paris-Saclay, CNRS, Laboratoire de Math\'ematiques d'Orsay (UMR 8628), 91405 Orsay Cedex, France (F.~Rousset)}
\email{frederic.rousset@universite-paris-saclay.fr}

\author{Katharina Schratz}
\address{LJLL (UMR 7598), Sorbonne Universit\'e, UPMC, 4 place Jussieu, 75005, Paris, France (K.~Schratz)}
\email{katharina.schratz@ljll.math.upmc.fr}
\thanks{KS gratefully acknowledges funding from the European Research Council (ERC) under the European Union’s Horizon 2020 research and innovation programme (grant agreement No.~850941).}

\begin{document}

\begin{abstract}
A filtered Lie splitting scheme is proposed for the time integration of the cubic nonlinear Schr\"{o}dinger equation on the two-dimensional torus $\mathbb{T}^2$. The scheme is analyzed in a framework of discrete Bourgain spaces, which allows us to consider initial data with low regularity; more precisely initial data in $H^s(\mathbb{T}^2)$ with $s>0$. In this way, the usual stability restriction to smooth Sobolev spaces with index $s>1$ is overcome. Rates of convergence of order $\tau^{s/2}$ in $L^2(\mathbb{T}^2)$ at this regularity level are proved. Numerical examples illustrate that these convergence results are sharp.
\end{abstract}

\maketitle

\section{Introduction}

We consider the cubic periodic nonlinear Schr\"{o}dinger equation (NLS)
\begin{equation}\label{n}
i\partial_t u=-\Delta u-\mu|u|^2u,\quad (t,x)\in\mathbb{R}\times\mathbb{T}^2
\end{equation}
on the two-dimensional torus $\mathbb{T}^2$. Our framework will allow us to study both the focusing case ($\mu=1) $ and the defocusing case ($\mu=-1$).

For the numerical integration of this equation, splitting methods are often the method of choice. For smooth initial data, their convergence behavior is nowadays well understood \cite{Faou12,Lubich08}. The standard Lie splitting method requires initial data in $H^2$ for proving first-order convergence in $L^2$. For less smooth initial data, one can prove error estimates in $L^2$ with the fractional rate of convergence $\tau^{s/2}$ for solutions in $H^s$ as long as $s>1$, see \cite{ESS16}. The restriction $s>1$ in two space dimensions is crucial due to the employed Sobolev embedding theorem
\begin{equation}\label{sob}
\Vert  v w \Vert_{L^2 } \leq C \Vert v  \Vert_{L^2}\Vert w \Vert_{H^{\sigma}},\quad \sigma>1.
\end{equation}
For $s<1$, new techniques have to be applied for estimating the nonlinear terms in the error recursion. In addition, schemes must be used that can handle different frequencies of the solution in different ways. For full space problems, a frequency filtered Lie splitting scheme was considered in \cite{Ignat11}. Using Strichartz estimates, convergence of order one in $L^2$ was shown for initial data in $H^2$. In our previous papers \cite{Ost, Ost1}, we considered Fourier integrators and a filtered Lie splitting scheme for the one-dimensional cubic nonlinear Schr\"{o}dinger equation on the torus. There, the analysis was carried out in a framework of discrete Bourgain spaces.

For low regularity integrators under Neumann boundary conditions we refer to \cite{BLW}; a low regularity integrator without loss of derivatives for cubic NLS in smooth Sobolev space was introduced in \cite{WY22} and very recently a modified Strang splitting was introduced in \cite{Wu} for the one dimensional cubic NLS equation which allows for improved error structure.

The filtered Lie splitting scheme
\begin{equation}\label{0}
u_{n+1}=\Psi^\tau(u_n)=e^{i\tau\Delta}\Pi_\tau(e^{\mu i\tau|\Pi_\tau u_n|^2}\Pi_\tau u_n),\quad u_0=\Pi_\tau u(0)
\end{equation}
is one of the simplest numerical schemes for NLS. Here the projection operator $\Pi_\tau$ for $\tau>0$ is defined by the Fourier multiplier:
\begin{equation}\label{proj}
\Pi_\tau=\bar{\Pi}_\tau=\chi\left(\dfrac{-i\nabla}{\tau^{-\frac{1}{2}}}\right),
\end{equation}
where $\chi$ is the characteristic function of the square $[-1,1]^2$. Note that the filter depends on the size of the time step $\tau$. This filter is crucial for the scheme to approximate well the frequency interactions of the continuous model. For the stability analysis at low regularity, which is carried out in this paper, the key trilinear estimate introduced in Theorem \ref{disckeythm} would not be true without this filter, or it would involve an additional loss of derivative with a filter which allows higher frequencies, see for example \cite{ORS1} for this phenomenon in the case of Strichartz estimates in the whole space.

The filtered Lie splitting scheme has shown good convergence properties for initial data with low regularity for the one-dimensional NLS \cite{Ost1}. In this paper, we will study the two-dimensional case and prove the following result:

\begin{theorem}\label{mainthm}
Let  $u_0\in H^{s_0}(\mathbb{T}^2),~s_0\in (0, 2]$ and $T>0$  such that there exists an exact solution $u$  of \eqref{n} with initial data $u_0$  such that  $u\in X^{s_0, b_0}(T)\subset\mathcal{C}([0,T],H^{s_0})$. Let $u_n$ be the numerical solution defined by the scheme \eqref{0}. Then, we have the following error estimate: there exist $\tau_0>0$ and $C_T>0$ such that for every time step size $\tau\in(0,\tau_0]$,
\begin{equation}\label{c}
\Vert u_n-u(t_n)\Vert_{L^2(\mathbb{T}^2)}\leq C_T\tau^{\frac{s_0}2},\quad 0\leq n\tau\leq T.
\end{equation}
\end{theorem}

The definition of the space $X^{s,b}(T)$ is given below in section~\ref{sectionanalysis}; the local Cauchy theory for rough data for 2D NLS is also recalled in this section. The key point is a subtle multilinear estimate in Bourgain spaces $X^{s,b}$. The analysis of the one-dimensional case was performed in \cite{Ost1} relying on the use of the discrete Bourgain spaces introduced in \cite{Ost}. There are some significant new difficulties in order to handle the 2D case in the analysis of the Cauchy problem. As the $L^4$ Strichartz estimate is not sufficient for the proof, we have to develop the appropriate multilinear estimates at the discrete level. We shall actually introduce a variant (see Lemma~\ref{mult2} and \eqref{rbis} for the discrete version) which allows us to analyze the error of the scheme in $L^2$ instead of $H^s$, $s>0$. Note that the result we obtain here, is slightly better than the one we obtained in \cite{Ost1} in the 1D case, where we got the convergence rate $\tau^{s \over 2}$ for every $s <s_{0}$ instead of $\tau^{\frac{s_0}2}$. Some of the improvements in the analysis  which are introduced here can be also used in the 1D case to get the sharp rate of convergence $\tau^{s_0\over 2}$.

In the defocusing case ($\mu=-1$), for  $s_{0} \geq 1$, the solutions are global. This is easily seen by using the conservation of the $L^2$ norm and the energy
$$
H(u)=  \int_{\mathbb{R}^2} {1 \over 2} | \nabla u|^2 + {1 \over 4}|u|^4  \, dx
$$
and hence  $T$ can be taken arbitrarily large in the above convergence result. This remains true if $s_{0}$ is close enough to $1$. For example, thanks to \cite{Pavlovic}, the solutions are global for $s_{0} \geq 2/3$. In the focusing case ($\mu=1$), this still holds if in addition the $L^2$ norm is small enough. However, blow-up solutions are known to exist. In this case, our convergence result is valid up to the blow-up.

To simplify the notations, we take from now on $\mu=-1$. We stress, however, that our analysis given below remains true for $\mu=1$ as long as the solution exists on $[0, T]$.

We note that \eqref{0} can be seen as a classical Lie splitting discretization of the projected equation
\begin{equation}\label{p}
i\partial_t u^\tau=-\Delta u^\tau+\Pi_\tau(|\Pi_\tau u^\tau|^2\Pi_\tau u^\tau), \quad u^\tau(0)=\Pi_\tau u(0).
\end{equation}
This property will be used in our proofs.

\subsection*{Outline of the paper.} The  paper is organized as follows. In section~\ref{sectionanalysis}, we briefly recall the main steps of the analysis of the Cauchy problem for \eqref{n}, and we prove an estimate on the difference between the exact solutions of \eqref{n} and \eqref{p}. In section~\ref{sectiondiscbourg}, we give the main properties of the discrete Bourgain spaces  $X^{s,b}_{\tau}$. The required technical Strichartz and multilinear estimates will be proven in section~\ref{sectionproof}. In section~\ref{sectioncontidisc}, we estimate the restriction of the continuous solution $u^{\tau}$ onto the grid in the framework of discrete Bourgain spaces. In section~\ref{sectionlocal}, we analyze the local error of the scheme \eqref{0}, and we give global error estimates in section~\ref{sectionglobal}. Finally we prove our main result, Theorem~\ref{mainthm}, in section~\ref{sectionmainthm}. Numerical examples, which are given at the end of the paper, illustrate our convergence result.

\subsection*{Notations.} The estimate $A\lesssim B$ means $A\leq CB$, where $C$ is a generic constant; in particular, $C$ is independent of the time step size $\tau \in (0, 1]$. The symbol $\lesssim_\gamma$ emphasizes that the constant $C$ depends on $\gamma$. Moreover, $A\sim B$ means that $A\lesssim B\lesssim A$.

We further denote $z \cdot w = \text{Re\hspace{1pt}} (z \,\overline{w})$  for $z, w \in \mathbb{C}$. This is the real scalar product in $\mathbb{C}$.
We also use the notation $\langle \,\cdot\, \rangle = ( 1 + | \cdot |^2)^{1 \over 2}$. Finally, for sequences $(u_{n})_{n \in \mathbb{Z}} \in  X^\mathbb{Z}$ in a Banach space $X$ with norm $\|\cdot \|_{X}$, we employ the usual norms
$$
\|u_{n}\|_{l^p_{\tau}X}=\left( \tau \sum_{n} \|u_{n}\|_{X}^p \right)^{1 \over p}, \quad
\|u_{n}\|_{l^\infty_{\tau}X} = \sup_{n \in \mathbb{Z}} \|u_{n}\|_{X}.
$$

\section{Analysis of the exact solution}\label{sectionanalysis}

In this section, we will discuss the Cauchy problem for \eqref{n} at low regularity. The use of Bourgain spaces and some subtle multilinear estimates is crucial for that. We shall then also discuss the Cauchy problem for the projected equation \eqref{p} and the error estimate between the solution $u$ of~\eqref{n} and the solution $u^\tau$ of \eqref{p}.

Let us first define Bourgain spaces. For a function $u(t,x)$ on $\mathbb{R}\times\mathbb{T}^2$, $\tilde{u}(\sigma,k)$ stands for its time-space Fourier transform, i.e.
$$
\tilde{u}(\sigma,k)=\int_{\mathbb{R}\times\mathbb{T}^2}u(t,x)e^{-i\sigma t-i\langle k,x\rangle}dx dt,
$$
where $\langle\cdot,\cdot\rangle$ denotes the inner product in $\mathbb R^2$. The inverse transform is given by
$$
u(t,x)=\frac{1}{4\pi^2}\sum\limits_{k\in\mathbb{Z}^2}\hat{u}_k(t)e^{i\langle k,x\rangle}
$$
with the Fourier coefficients $\hat{u}_k(t)=\frac{1}{2\pi}\int_{\mathbb{R}}\tilde{u}(\sigma,k)e^{i\sigma t}d\sigma$.

This way, we can define the Bourgain space $X^{s, b}= X^{s,b}(\mathbb{R}\times\mathbb{T}^2)$ consisting of functions with norm
$$
\Vert u\Vert_{X^{s,b}}=\Vert\langle k\rangle^s\langle\sigma + |k|^2\rangle^b\tilde{u}(\sigma,k)\Vert_{L^2 l^2}.
$$

We will also define a localized version of this space $ X^{s,b}(I)$ where $I\subset\mathbb{R}$ by
\begin{equation}\label{bourgainloc}
\Vert u\Vert_{X^{s,b}(I)}=\inf\{\Vert v\Vert_{X^{s,b}},~v|_I=u\},
\end{equation}
and we write $X^{s,b}(T)=X^{s,b}([0, T])$ for short.

We shall first recall some well-known properties of Bourgain spaces.

\begin{lemma}\label{contiprop}
For $\eta\in\mathcal{C}_c^\infty(\mathbb{R})$, we have that
\begin{align}
\Vert\eta(t)e^{it\Delta}f\Vert_{X^{s,b}}&\lesssim_{\eta,b}\Vert f\Vert_{H^s(\mathbb{T}^2)},\quad s\in\mathbb{R},~b\in\mathbb{R},~f\in H^s(\mathbb{T}^2);\nonumber\\
\Vert\eta(t)u\Vert_{X^{s,b}}&\lesssim_{\eta,b}\Vert u\Vert_{X^{s,b}},\quad s\in\mathbb{R},~b\in\mathbb{R};\nonumber\\
\Vert\eta\bigl(\tfrac{t}T\bigr)u\Vert_{X^{s,b^\prime}}&\lesssim_{\eta,b,b^\prime}T^{b-b^\prime}\Vert u\Vert_{X^{s,b}},\quad s\in\mathbb{R},~-\tfrac12<b^\prime<b<\tfrac12,~T\in(0,1];\nonumber\\
\left\Vert\eta(t)\int_{-\infty}^te^{i(t-t')\Delta}F(t')dt'\right\Vert_{X^{s,b}}&\lesssim_{\eta,b}\Vert F\Vert_{X^{s,b-1}},\quad s\in\mathbb{R},~b>\tfrac12;\nonumber\\
\label{d}\Vert u\Vert_{L^\infty H^s}&\lesssim_b\Vert u\Vert_{X^{s,b}},\quad s\in\mathbb{R},~b>\tfrac{1}{2}.
\end{align}
\end{lemma}

These estimates are proven, for example, in \cite[section~2.6]{Tao}.

A consequence of the above estimates and the definition of the local spaces \eqref{bourgainloc} is that  for every $b \in (\frac12, 1]$, $b'\in (0, \frac12)$ with $b+b'<1$ and $s \in \mathbb{R}$, we have uniformly for $0 <T \leq 1$,
\begin{equation}\label{bourgainloc1}
\Vert e^{it\Delta}f\Vert_{X^{s,b}(T)} \lesssim    \|f\|_{H^s}, \qquad   \left\Vert \int_{0}^t e^{i(t-t') \Delta} F(t') \, dt' \right\Vert_{X^{s,b}(T)} \lesssim  T^{1 - b- b'} \| F \|_{X^{s,-b'}(T)}.
\end{equation}

Besides, we need the following crucial  trilinear estimate for the analysis of \eqref{n}.

\begin{proposition}\label{propmult}
For any $s>0,~b>\max(\frac{1}{4},\frac{1}{2}-\frac{1}{4}s)$, and $u,v,w\in X^{s,b}$, we have that
\begin{equation}\label{s}
\Vert u \overline{v}w\Vert_{X^{s,-b}}\lesssim\Vert u\Vert_{X^{s,b}}\Vert v\Vert_{X^{s,b}}\Vert w\Vert_{X^{s,b}}.
\end{equation}
\end{proposition}

The main ingredients for the proof of this crucial estimate will be recalled later when we study the discrete version of this inequality.
Note that by using again the definition \eqref{bourgainloc} of the local spaces, we also have for $s, \, b$ as above and  for every $T>0$
\begin{equation}
\label{theomult1}\Vert u\overline{v} w\Vert_{X^{s,-b}(T)}\lesssim\Vert u\Vert_{X^{s,b}(T)}\Vert v\Vert_{X^{s,b}(T)}\Vert w\Vert_{X^{s,b}(T)}.
\end{equation}
The above estimate is uniform in $T$ for $T \in (0, 1]$.

We shall now establish the existence and uniqueness of the exact solution to \eqref{n}.

\begin{theorem}\label{theoexist}
Let $s_0>0$ and $u_0\in H^{s_{0}}(\mathbb{T}^2).$ Then, there exists  $T^*\in (0, +\infty]$ and a unique maximal solution $u\in  X^{s_0,b_{0}}(T) \subset \mathcal{C}([0,T],H^{s_0}(\mathbb{T}^2))$ for every $T \in (0, T^*)$ satisfying \eqref{n} for any $b_{0}\in\left(\frac{1}{2},\min(\frac{1}{2}+\frac{1}{4}s_0,\frac{3}{4})\right)$.
\end{theorem}

\begin{proof}
Let $s_0>0,~b_{0}\in\left(\frac{1}{2},\min(\frac{1}{2}+\frac{1}{4}s_0,\frac{3}{4})\right)$, and consider the Duhamel form of \eqref{n}
$$
v= F(v),
$$
where
\begin{equation}\label{x}
F(v)(t)=e^{it\Delta}u_0-i \int_0^te^{i(t-\vartheta)\Delta}\Big(|v(\vartheta)|^2v(\vartheta)\Big)d\vartheta.
\end{equation}
Thanks to \eqref{bourgainloc1} we have uniformly for every $T_{1} \in (0, 1]$
$$
\|F(v)\|_{X^{s_0, b_0}(T_{1})} \lesssim \|u_0\|_{H^{s_0}} + T_1^{\varepsilon_0} \||v|^2v\|_{ X^{s_0,-b'}(T_1)},
$$
where $b'$ is chosen such that $\max(\frac14,\frac12-\frac14s_0)<b^\prime<1-b_0<\frac12$ (so that $b^\prime<\frac12<b_0$) and $\varepsilon_0=1-b_0-b'>0$. By using the local version \eqref{theomult1} of Proposition~\ref{propmult}, we thus deduce that
$$
\Vert F(v)\Vert_{X^{s_0,b_{0}}(T_{1})}\leq C\Vert u_0\Vert_{H^{s_0}}+CT_{1}^{\varepsilon_0}\Vert v\Vert_{X^{s_0,b_{0}}(T_{1})}^3.
$$
From the same arguments, we can get
\begin{equation}\label{contraction}
\Vert F(v_1)-F(v_2)\Vert_{X^{s_0,b_{0}}(T_{1})}\leq 3CT_{1}^{\varepsilon_0}R^2\Vert v_1-v_2\Vert_{X^{s_0,b_{0}}(T_{1})}\end{equation}
for all $v_1,v_2$ satisfying $\Vert v_i\Vert_{X^{s_0,b_{0}}(T_{1})}\leq R$. Note that in the above estimates  $C$ is independent of $T_{1}\in (0,1]$.
By fixing $R>  C\Vert u_0\Vert_{H^{s_0}}$, we get by the Banach fixed point theorem the existence of  a fixed point of $F$ (and thus of a solution of the PDE) in the ball of radius $R$ of $X^{s_{0}, b_{0}}(T_{1})$ for $T_{1}$ sufficiently small. The fact that $X^{s_{0}, b_{0}}(T_{1})$ is contained in $\mathcal{C}([0, T_1], H^{s_0})$ is a standard property of Bourgain spaces for $b_0>1/2$ (this is actually obtained while proving \eqref{d}).

From the same arguments we can also get that for every $T>0$ there is  uniqueness of solutions in $X^{s_{0}, b_{0}}(T)$ for $s_{0}$ and $b_{0}$ as above. Indeed, if $v_{1}$ and $v_{2}$ are two solutions in  $X^{s_{0}, b_{0}}(T)$, we can set
$$
R_{T}= \|v_{1}\|_{X^{s_{0}, b_{0}}(T)} +  \|v_{2}\|_{X^{s_{0}, b_{0}}(T)}.
$$
We can then cover $[0,T]$ by a finite number of intervals of the form $[nT_{1}, nT_{1}+T_{1}] \cap [0,T]$, where $T_{1}>0$ is chosen such that $CT_{1}^{\varepsilon_0} R_{T}^2\leq 1/2$. We then obtain easily from a shifted version of the estimate \eqref{contraction}  that if $v_{1}$ and $v_{2}$ coincide at $nT_{1}$ then they coincide on $[nT_{1}, nT_{1}+T_1]$ and thus get that $v_{1}=v_{2}$ on $[0, T]$.

From the existence and  uniqueness, we get the  maximal solution in a standard way.
\end{proof}

The next step will be to study the projected equation \eqref{p}.

\begin{theorem}\label{theou-utau}
Let $u_{0}\in H^{s_{0}}$ for some $s_{0}>0$ and  $u  \in X^{s_{0},b_{0}}(T)$ for any  $T\in (0, T^*)$ the solution of \eqref{n} given by Theorem~\ref{theoexist}. Then there exists $\tau_{0}>0$ such that for every $\tau \in (0, \tau_{0}]$, there exists a unique solution $u^\tau$ of \eqref{p} which is also defined on $[0, T]$ and such that $u^\tau \in X^{s_{0}, b_{0}}(T)$. Moreover, we also have that
\begin{equation} \label{diffutau}
\|u- u^\tau \|_{X^{0, b_0}(T)}  \lesssim C_{T} \tau^{s_0\over 2}, \qquad \|u^\tau \|_{X^{s_0, b_0}(T)} \leq C_T
\end{equation}
for some $C_T$ independent of $\tau$ for $\tau \in (0,\tau_0].$
\end{theorem}

\begin{proof}
The existence part for $u^\tau$ is very easy. Indeed from the smoothing effect of the projection $\Pi_{\tau}$, \eqref{p} can be seen as an ordinary differential equation on $L^2$ and thus by  the conservation of the $L^2$ norm,  we get that \eqref{p} is globally  well-posed in $L^2$ and the obtained solution is actually $C^\infty$. The difficult part is then to prove the uniformity in $\tau$ of the estimates \eqref{diffutau} on $[0, T]$. Let us set $v= u- u^\tau$, we shall first prove that for $\tau$ sufficiently small, $v$ is uniformly bounded in $X^{s_0, b_0}(T)$. Since $ u\in X^{s_0,b_0}(T)$ by assumption, this will give the second part of \eqref{diffutau}. By easy algebraic computations, we get that $v$ solves
$$
i \partial_{t}v + \Delta v = L(u,v) + Q(u, v) + C(u, v) + S(u), \quad v|_{t=0}= (I-\Pi_{\tau})u_{0},
$$
where
\begin{align*}
L(u, v)&=  \Pi_{\tau}\left( |\Pi_{\tau}u|^2 \Pi_{\tau}v  + 2 \Pi_{\tau}u \cdot \Pi_{\tau}v \Pi_{\tau}u\right), \\
Q(u,v)&= -\Pi_{\tau}\left( |\Pi_{\tau}v|^2 \Pi_{\tau}u + 2 \Pi_{\tau}u \cdot \Pi_{\tau}v \Pi_{\tau}v \right), \\
C(u, v) &= \Pi_{\tau} \left( | \Pi_{\tau} v|^2 \Pi_{\tau}v \right), \\
S(u) &=  (I-\Pi_{\tau})( |u|^2 u) + \Pi_{\tau}\Bigl( |\Pi_{\tau}u |^2 (I-\Pi_{\tau})u + 2 \Pi_{\tau} u\cdot  (I-\Pi_{\tau}) u  \Pi_{\tau}u \\[-2mm]
 & \quad{}+ \Pi_{\tau}u | (I-\Pi_{\tau})u|^2 + 2 \Pi_{\tau}u \cdot (I-\Pi_{\tau})u (I-\Pi_{\tau}) u + |(I-\Pi_{\tau}) u |^2 (I-\Pi_{\tau})u\Bigr).
\end{align*}
Duhamel's formula thus reads
\begin{equation} \label{Duhamelbis}
v(t) = e^{it\Delta} (I-\Pi_{\tau})u_{0} - i \int_{0}^t  e^{i(t-t') \Delta}\left( L(u,v) (t') + Q(u,v) (t') + C(u,v)(t') + S(u)(t') \right) dt'.
\end{equation}
Let us choose $b'$ as in the proof of Theorem~\ref{theoexist}. We observe that $\Pi_{\tau}$ is continuous on $X^{s,b}$ for every $s$, $b$, and we have that $\|\Pi_\tau u\|_{X^{s,b}(T)} \leq \|u\|_{X^{s,b}(T)}$ uniformly for $\tau\in (0,1]$, $T>0$ and $u\in X^{s,b}(T)$. Moreover, $\Pi_{\tau}u$ converges to $u$ for any $u$ in $X^{s,b}(T)$ for every $s$, $b$. Thus we obtain that for every $\varepsilon>0$, there exists $\tau_0>0$ (depending on $T$ and $\varepsilon$) such that for every $\tau\in(0,\tau_0]$, we have
\begin{equation}\label{sourcesmall1}
\| e^{it\Delta} (I-\Pi_{\tau})u_{0}\|_{X^{s_{0}, b_{0}}(T)} \lesssim \| (I-\Pi_{\tau})u_{0}\|_{H^{s_{0}}} \leq \varepsilon, \qquad  \| (I-\Pi_{\tau})u\|_{X^{s_{0}, b_{0}}(T)} \leq \varepsilon.
\end{equation}
Indeed, since $I-\Pi_{\tau}$ projects on spatial frequencies higher than $\tau^{-\frac12},$ this follows from the definitions of the norms and the dominated convergence theorem. Moreover, by using Proposition~\ref{propmult} and by choosing $\tau_{0}$ maybe smaller, we can also achieve
\begin{multline}\label{sourcesmall2}
\left\| \int_{0}^t  e^{i(t-t') \Delta}  S(u)(t') \, dt' \right\|_{X^{s_{0}, b_{0}}(T)} \lesssim_{T} \|(I-\Pi_{\tau}) u \|_{X^{s_{0}, b_{0}}(T)} \|u \|_{X^{s_{0}, b_{0}}(T)}^2  \\+\| (I- \Pi_{\tau}) ( |u|^2u ) \|_{X^{s_{0}, -b'}(T)}  \leq \varepsilon.
\end{multline}
Let us  take $T_{1} \in (0, \min (1, T)]$ (to be chosen sufficiently small later), and set
$$
R_{T}= \|u\|_{X^{s_{0}, b_{0}}(T)}.
$$
From \eqref{Duhamelbis}, by using again \eqref{bourgainloc1}, \eqref{theomult1}, \eqref{sourcesmall1}, and \eqref{sourcesmall2}, we obtain the estimate (which is uniform in $T_1\in(0,\min(1,T)]$ due to the uniformity of \eqref{bourgainloc1})
$$
\|v\|_{X^{s_0, b_0}(T_1)} \leq  2 \varepsilon + T_1^{\varepsilon_0} R_T^2 \|v\|_{X^{s_0,b_0}(T_1)} + R_T\|v\|_{X^{s_0, b_0}(T_1)}^2 + C_T\|v\|_{X^{s_0, b_0}(T_1)}^3
$$
for some uncritical number $C_T>0$. Consequently, by choosing at first $T_1$ sufficiently small so that $T_{1}^{\varepsilon_0} R_T^2 \leq 1/2$ and then $\varepsilon$ sufficiently small (and so $\tau_0$ sufficiently small), we can deduce that
$$
\|v\|_{X^{s_{0}, b_{0}}(T_{1})} \lesssim  8\varepsilon.
$$
In fact, let us take $\varepsilon$ sufficiently small so that
$$
R_{T} \varepsilon^{1 \over 2 } + C_{T} \varepsilon  \leq {1 \over 8},
$$
we can easily get the existence of a fixed point for \eqref{Duhamelbis}
in the ball  $\|v\|_{X^{s_{0}, b_{0}}(T_{1})} \leq \varepsilon^{1 \over 2}$
for $\varepsilon$ sufficiently small from the above estimates.

Since the choice of $T_{1}$ depends only on $R_{T}$, we can iterate the argument on $[T_{1}, 2T_{1}]$ and so on to cover $[0, T]$ and get that $\|v\|_{X^{s_{0}, b_{0}}(T)}$ is bounded on $[0, T]$. Since $u \in X^{s_{0}, b_{0}}(T)$, this yields the second part of \eqref{diffutau}. It remains to prove the first part of \eqref{diffutau}. The proof will rely on a non-standard variant of \eqref{s}.

\begin{lemma}\label{mult2}
For any $s>0,~b>\max(\frac{1}{4},\frac{1}{2}-\frac{1}{4}s)$, and $\sigma$ any permutation of $\{1, 2, 3 \}$,
if  $u_{\sigma(1)},u_{\sigma(2)}\in X^{s,b}$, $u_{\sigma(3)} \in X^{0,b}$,  we have that
\begin{equation}\label{eq:mult2}
\Vert u_{1} \overline{u_{2}} u_{3}\Vert_{X^{0,-b}}\lesssim \Vert u_{\sigma(1)}\Vert_{X^{s,b}}\Vert u_{\sigma(2)}\Vert_{X^{s,b}}\Vert u_{\sigma(3)}\Vert_{X^{0,b}}.
\end{equation}
\end{lemma}

Note that this means that when we perform the estimate with zero space regularity, we can avoid to lose derivatives on one of the three functions (which can by any of them). We shall give below the details of the proof of the discrete counterparts of \eqref{s} and the above estimate. The same techniques can be used to deduce \eqref{eq:mult2}. From the above estimate and the definition of the local Bourgain spaces, we also deduce the local version: for every $T>0$
\begin{equation}\label{theomult2}
\Vert u_1 \overline{u_2}u_{3}\Vert_{X^{0,-b}(T)}\lesssim\Vert u_{\sigma(1)}\Vert_{X^{s,b}(T)}\Vert u_{\sigma(2)}\Vert_{X^{s,b}(T)}\Vert u_{\sigma(3)}\Vert_{X^{0,b}(T)}.
\end{equation}
Again, the estimate is uniform in $T$ for $T \in (0, 1]$.

Now let us observe that again by the properties of $\Pi_{\tau}$, we have
\begin{equation} \label{sourcesmall1bis}
\| e^{it\Delta} (I-\Pi_{\tau})u_{0}\|_{X^{0, b_{0}}(T)} \lesssim  \| (I-\Pi_{\tau})u_{0}\|_{L^2} \lesssim \tau^{s_{0}\over 2}, \qquad  \| (I-\Pi_{\tau})u\|_{X^{0, b_{0}}(T)} \lesssim \tau^{{s_{0}\over 2}}
\end{equation}
since $u \in  X^{s_{0}, b_{0}}(T)$ and $u_{0}\in H^{s_{0}}$ where again, the involved constant depends on $T$.

By taking $b'$ as before, we get from \eqref{theomult2} that
\begin{equation}\label{sourcesmall2bis}
\begin{aligned}
\Big\| \int_{0}^t  e^{i(t-t') \Delta} & S(u)(t') \, dt'
\Big\|_{X^{0, b_{0}}(T)}\\
&\lesssim_{T} \|(I-\Pi_{\tau}) u \|_{X^{0, b_{0}}(T)} \|u \|_{X^{s_{0}, b_{0}}(T)}^2  +\| (I- \Pi_{\tau}) ( |u|^2u ) \|_{X^{0, -b'}(T)}\\
&\lesssim_{T} \tau^{{s_{0}\over 2}}  \| u  \|^3_{X^{s_{0}, b_0}(T)} + \tau^{{s_{0}\over 2}} \| |u|^2u  \|_{X^{s_{0}, -b'}(T)}\lesssim_{T} \tau^{{s_{0}\over 2}}
\end{aligned}
\end{equation}
since $|u|^2u \in X^{s_{0}, -b'}(T)$ thanks to \eqref{theomult1}.

Now, let us take
$$
\widetilde R_{T}= \|u \|_{X^{s_{0}, b_{0}}(T)} + \|v\|_{X^{s_{0}, b_{0}}(T)}
$$
which is well defined since we have already proven that  $\|v\|_{X^{s_{0}, b_{0}}(T)}$ is bounded. By using \eqref{bourgainloc1}, \eqref{Duhamelbis}, \eqref{sourcesmall1bis}, \eqref{sourcesmall2bis}, and \eqref{theomult2}, we now get for every $T_{1} \leq \min (1, T)$,
\begin{align*}
\|v\|_{X^{0, b_{0}}(T_{1})} &\lesssim_{T} \tau^{s_{0} \over 2} + T_{1}^{\varepsilon_{0}} \|u\|_{X^{s_{0}, b_{0}}(T)}^2   \|v\|_{X^{0, b_{0}}(T_{1})}\\
&\qquad{}+T_{1}^{\varepsilon_{0}} \|u\|_{X^{s_{0}, b_{0}}(T)} \| v \|_{X^{s_{0}, b_{0}}(T)}    \|v\|_{X^{0, b_{0}}(T_{1})}  + T_{1}^{\varepsilon_{0}} \| v \|_{X^{s_{0}, b_{0}}(T)}^2  \|v\|_{X^{0, b_{0}}(T_{1})} .
\end{align*}
This yields
$$
\|v\|_{X^{0, b_{0}}(T_{1})} \lesssim_{T} \tau^{s_{0} \over 2}  + T_{1}^{\varepsilon_0} { \widetilde R_T^2 }  \|v\|_{X^{0, b_{0}}(T_{1})}
$$
and hence, by taking $T_{1}$ sufficiently small such that  $ T_{1}^{\varepsilon_{0}} \widetilde R_{T}^2  \leq \frac12$ (thus $T_{1}$ depends only on $T$), we get that
$$
\|v\|_{X^{0, b_{0}}(T_{1})} \lesssim_{T} \tau^{s_{0} \over 2}.
$$
We can then  again  iterate the argument on $[T_{1}, 2T_{1}]$ and so on to cover $[0, T]$ and get  that
$$
\|v\|_{X^{0, b_{0}}(T)} \lesssim_{T} \tau^{s_{0} \over 2}.
$$
This ends the proof of Theorem~\ref{theou-utau}.
\end{proof}

The final result of this section is an estimate of $u^\tau$ in Bourgain spaces with larger index $b$. This will be useful to prove the boundedness of $u^\tau$ in the discrete spaces.

\begin{corollary}\label{corhigherb}
Under the assumptions of Theorem~\ref{theou-utau}, we have for every $b_{1}>\frac12$ such that $b_{1}\leq b_{0}$ and for every  $\sigma >\frac13$ the estimate
\begin{equation}\label{utaumieux}
\|u^\tau \|_{X^{s_{0},  1-b_{0}+ b_{1}}(T)} \lesssim_{T, \sigma} \tau^{-\frac{3\sigma}2}.
\end{equation}
\end{corollary}

Note that $\sigma$ will be chosen sufficiently close to $\frac13$ later.

\begin{proof}
Let us recall that $u^\tau\in X^{s_0, b_0}(T)$ solves
$$
u^\tau(t)=e^{it\Delta}u_0+\int_0^te^{i(t-\vartheta)\Delta} F(u^\tau(\vartheta))d\vartheta, \quad F(v)=-i\Pi_\tau (|\Pi_\tau v|^2\Pi_\tau v).
$$
Let us take $b_1>\frac12$. We shall estimate $\|u^\tau \|_{X^{s_0, 1-b_0+ b_1}}$ which is nontrivial when   $1-b_0+b_1>b_0$. We first write that
$$
\|u^\tau\|_{X^{s_0, 1-b_0+b_1}(T)}\lesssim\|u_0\|_{H^{s_0}}+\Vert F(u^\tau)\Vert_{X^{s_0,1-b_0+b_1-1}(T)}
\lesssim\|u_0\|_{H^{s_0}}+\Vert\langle\partial_x\rangle^{s_0}F(u^{\tau})\Vert_{L^2([0,T]\times\mathbb{T}^2)}
$$
since $1-b_0+b_1-1\leq0$. From the generalized Leibniz rule (see \cite{Mus}), we thus get
$$
\Vert F(u^{\tau})\Vert_{X^{s_{0},1-b_{0}+b_{1}-1}(T)}\lesssim \Vert\langle\partial_x\rangle^{s_0}\Pi_\tau u^{\tau}\Vert_{L^6([0, T] \times \mathbb{T}^2)}^3.
$$
By the $L^6$ Strichartz estimate obtained in \cite{Bourg}, for any $\sigma>\frac13$ and $b_2>\frac12$, we have
$$
\Vert\langle\partial_x\rangle^{s_0}\Pi_\tau u^\tau\Vert_{L^6([0,T] \times \mathbb{T}^2)}\lesssim\Vert\Pi_\tau u^{\tau}\Vert_{X^{\sigma+s_0, b_2}(T)}.
$$
We thus get in particular, since $b_0>\frac12$,  that
$$
\|u^\tau\|_{X^{s_0,1-b_0+ b_1}(T)} \lesssim \|u_0\|_{H^{s_0}} + \|\Pi_\tau u^\tau\|_{X^{s_0+\sigma,b_0}(T)}^3
\lesssim \|u_0\|_{H^{s_0}} + \tau^{-\frac{3 \sigma} 2} \|u^\tau\|^3_{X^{s_0, b_0}(T)}
$$
by using again that $\Pi_{\tau}$ is projecting on spatial frequencies smaller than $\tau^{-\frac12}$. Consequently, by using \eqref{diffutau}, we get
$$
\|u^\tau \|_{X^{s_0, 1-b_0+ b_1}(T)} \lesssim_{T,\sigma}\tau^{-\frac{3\sigma}2},
$$
where $\sigma$ can be chosen arbitrarily close to $\frac13$. This ends the proof.
\end{proof}

\begin{remark}\label{remextension}
In the following, thanks to the definition of the local Bourgain spaces, we shall  still denote by $u^\tau$ an extension of the solution $u^\tau$ of \eqref{p} on $[0, T]$ such that
\begin{equation}\label{bmieux}
\|u^\tau \|_{X^{s_{0}, b_{0}}} \leq 2  \|u^\tau \|_{X^{s_{0}, b_{0}}(T)} \lesssim C_{T}, \qquad \|u^\tau \|_{X^{s_{0}, 1-b_{0}+ b_{1}}} \leq 2  \|u^\tau \|_{X^{s_{0}, 1-b_{0}+ b_{1}}(T)} \lesssim  C_{T, \sigma} \tau^{-\frac{3\sigma}2}
\end{equation}
with $\sigma$ close to $\frac13$ to be chosen. Note that $u^\tau$ is now defined globally, but it is a solution of \eqref{p} only on $[0, T]$.
\end{remark}

\section{Discrete Bourgain spaces}\label{sectiondiscbourg}

We now define the discrete Bourgain spaces, see also \cite{Ost}. We first take $(u_n(x))_n$ to be a sequence of functions on the torus $\mathbb{T}^2$, with its ``time-space" Fourier transform
$$
\widetilde{u_n}(\sigma,k)=\tau\sum\limits_{m\in\mathbb{Z}}\widehat{u_m}(k)e^{im\tau\sigma},
$$
where
$$
\widehat{u_m}(k)=\dfrac{1}{4\pi^2}\int_{\mathbb{T}^2}u_m(x)e^{-i\langle k,x\rangle}dx.
$$
In this framework, $\widetilde{u_{n}}$ is a $ 2 \pi/\tau$ periodic function in $\sigma$ and Parseval's identity reads
\begin{equation}\label{parseval}
\| \widetilde{u_{n}}\|_{L^2l^2}= \|u_{n}\|_{l^2_{\tau}L^2},
\end{equation}
where we use the shorthands
$$
\| \widetilde{u_n}\|_{L^2l^2}^2 = \int_{-{\pi \over \tau}}^{\pi\over\tau} \sum_{k \in\mathbb{Z}^2}
|\widetilde{u_n}(\sigma, k)|^2 \, d \sigma, \qquad
\|u_n\|_{l^2_{\tau}L^2}^2 = \tau \sum_{m \in \mathbb{Z}} \int_{\mathbb{T}^2} |u_m(x)|^2 \, dx.
$$

Then the discrete Bourgain space $X^{s,b}_\tau$ can be defined with the norm
\begin{equation}\label{5}
\Vert u_n\Vert_{X^{s,b}_\tau}=\Vert\langle k\rangle^s\langle d_\tau(\sigma-|k|^2)\rangle^b\widetilde{u_n}(\sigma,k)\Vert_{L^2l^2},
\end{equation}
where $d_\tau(\sigma)=\frac{e^{i\tau\sigma}-1}{\tau}$. Note that for any fixed $(u_n)_n$, the norm is an increasing function for both $s$ and $b$.

We can also define the localized discrete Bourgain spaces, $X_\tau^{s,b}(I)$, with
$$
\Vert u_n\Vert_{X_\tau^{s,b}(I)}=\inf\{\Vert v_n\Vert_{X_\tau^{s,b}}\,|\,v_n =u_n, \, n \tau \in I\}.
$$

We directly get from the definition \eqref{5} the elementary properties that for any $s\geq s^\prime$ and $b\geq b^\prime$, we have
\begin{align}
\label{a}\sup\limits_{\delta\in[-4,4]}\Vert e^{i\tau\delta\Delta} u_n\Vert_{X^{s,b}_\tau}&\lesssim\Vert u_n\Vert_{X^{s,b}_\tau},\\
\label{b}\Vert\Pi_\tau u_n\Vert_{X^{s,b}_\tau}&\lesssim\tau^{b^\prime-b}\Vert\Pi_\tau u_n\Vert_{X^{s,b^\prime}_\tau},\\
\label{t}\Vert\Pi_\tau u_n\Vert_{X^{s,b}_\tau}&\lesssim\tau^{\frac{s^\prime-s}{2}}\Vert\Pi_\tau u_n\Vert_{X^{s^\prime,b}_\tau}.
\end{align}
Another equivalent norm on  the discrete Bourgain spaces is given by the following characterization.

\begin{lemma}
For $(u_n)_n\in X^{s,b}_\tau$, let
$$
\Vert e^{-in\tau\Delta}u_n\Vert_{H_\tau^bH^s}=\Vert\langle D_\tau\rangle^b\langle\Delta\rangle^{\frac{s}{2}}(e^{-in\tau\Delta}u_n)\Vert_{l_\tau^2L^2},
$$
where $(D_\tau(u_n))_n=\big(\frac{u_{n-1}-u_n}{\tau}\big)_n.$ Then,
$$
\Vert u_n\Vert_{X^{s,b}_\tau}\sim\Vert e^{-in\tau\Delta}u_n\Vert_{H_\tau^bH^s},\quad u_n\in X^{s,b}_\tau,
$$
i.e., it is an equivalent norm in $X^{s,b}_\tau$.
\end{lemma}

\begin{proof}
First, by a change of variables, we have
$$
\Vert u_n\Vert_{X^{s,b}_\tau}^2=\sum\limits_{k\in\mathbb{Z}^2}\int_{-\frac{\pi}{\tau}}^{\frac{\pi}{\tau}}\langle d_\tau(\sigma)\rangle^{2b}\langle k\rangle^{2s}|\widetilde{u_n}(\sigma+|k|^2,k)|^2d\sigma.
$$
Then by setting $f_n(x)=e^{-in\tau\Delta}u_n(x)$, from the definition of the Fourier transform, we have
$$
\widetilde{f_n}(\sigma,k)=\tau\sum\limits_{m\in\mathbb{Z}}\widehat{u_m}(k)e^{im\tau(\sigma+|k|^2)},
$$
which implies
$$
\widetilde{f_n}(\sigma,k)=\widetilde{u_n}(\sigma+|k|^2,k).
$$
Moreover, by the definition of the Fourier transform, we have
$$
\widetilde{D_\tau(u_n)}(\sigma,k)=d_\tau(\sigma)\widetilde{u_n}(\sigma,k).
$$
Thus by Plancherel,
$$
\Vert u_n\Vert_{X^{s,b}_\tau}\sim\Vert\langle k\rangle^s\langle d_\tau(\sigma)\rangle^b\widetilde{f_n}\Vert_{L^2((-\tfrac\pi\tau,\tfrac\pi\tau)\times\mathbb{Z}^2)}\sim\Vert\langle D_\tau\rangle^b\langle\Delta\rangle^{\frac{s}{2}}f_n\Vert_{l_\tau^2L^2},
$$
which is the desired result.
\end{proof}

Next, we will give the counterparts of Lemma~\ref{contiprop} and Proposition~\ref{propmult}.
\begin{lemma}\label{discprop}
For $s\in\mathbb{R},~\eta\in\mathcal{C}_c^\infty(\mathbb{R})$ and $\tau\in(0,1]$, we have that
\begin{align}
\Vert\eta(n\tau)e^{in\tau\Delta}f\Vert_{X^{s,b}_\tau}&\lesssim_{\eta,b}\Vert f\Vert_{H^s(\mathbb{T}^2)},\quad b\in\mathbb{R},~f\in H^s,\nonumber\\
\Vert\eta(n\tau)u_n\Vert_{X^{s,b}_\tau}&\lesssim_{\eta,b}\Vert u_n\Vert_{X^{s,b}_\tau},\quad b\in\mathbb{R},~u_n\in X^{s,b}_\tau,\nonumber\\
\label{9}\Vert\eta(\tfrac{n\tau}{T})u_n\Vert_{X^{s,b^\prime}_\tau}&\lesssim_{\eta,b,b^\prime}T^{b-b^\prime}\Vert u_n\Vert_{X^{s,b}_\tau},\quad -\tfrac{1}{2}<b\leq b^\prime<\tfrac{1}{2},~0<T=N\tau\leq1,~N\geq1,\\
\label{y}\Vert u_n\Vert_{l_\tau^\infty H^s}&\lesssim_b\Vert u_n\Vert_{X^{s,b}_\tau},\quad b>\tfrac{1}{2},\\
\label{q}\Vert U_n\Vert_{X^{s,b}_\tau}&\lesssim_{\eta,b}\Vert u_n\Vert_{X^{s,b-1}_\tau},\quad b>\tfrac{1}{2},
\end{align}
where
$$U_n(x)=\tau\eta(n\tau)\sum\limits_{m=0}^ne^{i(n-m)\tau\Delta}u_m(x).$$
\end{lemma}

We stress all these estimates are uniform in $\tau$.

The proof of this lemma is also nearly the same as the one-dimensional case given in \cite[section~3]{Ost}. Therefore, we omit the details.

The key  nonlinear estimates for the analysis of the scheme are given in the following theorem.
\begin{theorem}\label{disckeythm}
For any $s>0,~b>\max(\frac{1}{4},\frac{1}{2}-\frac{1}{4}s),~b_1>\frac{1}{2}$, we have
\begin{align}
\label{z}\Vert\Pi_\tau u_n\Vert_{l^4_\tau L^4}&\lesssim\Vert u_n\Vert_{X_\tau^{s,b_1}},\\
\label{r}\Vert\Pi_\tau(\Pi_\tau u_{n,1} \overline{\Pi_\tau u_{n, 2}}\Pi_\tau u_{n,3})\Vert_{X_\tau^{s,-b}}&\lesssim\Vert u_{n,1}\Vert_{X_\tau^{s,b}}\Vert u_{n,2}\Vert_{X_\tau^{s,b}}\Vert u_{n,3}\Vert_{X_\tau^{s,b}}, \\
\label{rbis}\Vert\Pi_\tau(\Pi_\tau u_{n,1} \overline{\Pi_\tau u_{n, 2}}\Pi_\tau u_{n,3})\Vert_{X_\tau^{0,-b}}&\lesssim\Vert u_{n, \sigma(1)}\Vert_{X_\tau^{s,b}}\Vert u_{n,\sigma(2)}\Vert_{X_\tau^{s,b}}\Vert u_{n, \sigma(3)}\Vert_{X_\tau^{0,b}},
\end{align}
where $(u_n)_n$, $(u_{n,i})_n$, are functions in the corresponding spaces and $\sigma$ is any permutation of $\{1, 2, 3\}.$
\end{theorem}

We postpone the proof of Theorem~\ref{disckeythm} to section~\ref{sectionproof}. Note that \eqref{rbis} is the discrete counterpart of the estimate of Lemma~\ref{mult2}. The refinement in the proof of \eqref{r} that we shall use to get \eqref{rbis}, can be also easily performed at the continuous level to get Lemma~\ref{mult2}.

\section{Boundedness of the exact solution in discrete Bourgain spaces}\label{sectioncontidisc}

In this section, we shall prove the boundedness of $u^\tau(t_n)$ in the norm of $X^{s,b}_\tau$ for suitable $s$ and~$b$. Note that we denote here by $u^\tau$ the extension of the solution of \eqref{p} as defined in Remark~\ref{remextension}.

We first prove the following lemma.
\begin{lemma}\label{contidisc}
For any $s>0,~b<\frac{1}{2},~b^\prime>\frac{1}{2}$, and a given sequence of functions $(u_n(x))_{n\in\mathbb{Z}}$ with $u_n(x)=u(n\tau,x)$, we have
$$
\Vert u_n\Vert_{X_\tau^{s,b}}\lesssim\Vert u\Vert_{X^{s,b}}+\tau^{b^\prime}\Vert u\Vert_{X^{s,b+b^\prime}}.
$$
\end{lemma}

\begin{proof}
We adapt the proof in \cite[Lemma~3.4]{Rou}. We shall just prove the case $s=0$, the extension to general $s$ is straightforward. By setting $f=e^{-it\Delta}u$ and $f_n(x)=f(n\tau,x)$, it thus suffices to prove that
$$
\Vert f_n\Vert_{H_\tau^bL^2(\mathbb{Z}\times\mathbb{T}^2)}\lesssim\Vert f\Vert_{H^bL^2}+\tau^{b^\prime}\Vert f\Vert_{H^{b+b^\prime}L^2}.
$$
Since by definition,
$$
\widetilde{f_m}(\sigma,k)=\tau\sum\limits_{n\in\mathbb{Z}}\tilde{f}(n\tau,k)e^{in\tau\sigma},
$$
by Poisson's summation formula, we have
$$
\widetilde{f_n}(\sigma,k)=\sum\limits_{m\in\mathbb{Z}}\tilde{f}(\sigma+\tfrac{2\pi}{\tau}m,k).
$$
Therefore,
$$
\langle d_\tau(\sigma)\rangle^b\widetilde{f_n}(\sigma,k)=\sum\limits_{m\in\mathbb{Z}}\langle d_\tau(\sigma+\tfrac{2\pi}{\tau}m)\rangle^b\tilde{f}(\sigma+\tfrac{2\pi}{\tau}m,k)
$$
since $d_\tau$ is also a $2\pi/\tau$ periodic function. Since we always have that $|d_\tau(\sigma)|\lesssim\langle\sigma\rangle$ and $\tau\sigma\in[-\pi,\pi]$, by Cauchy--Schwarz,
\begin{align*}
\big|\langle d_\tau(\sigma)\rangle^b&\widetilde{f_n}(\sigma,k)\big|^2\\
&=|\langle d_\tau(\sigma)\rangle^b\tilde{f}(\sigma,k)|^2+\Big|\sum\limits_{m\neq0}\langle d_\tau(\sigma+\tfrac{2\pi}{\tau}m)\rangle^b\tilde{f}(\sigma+\tfrac{2\pi}{\tau}m,k)\Big|^2\\
&\lesssim|\langle\sigma\rangle^b\tilde{f}|^2+\Big|\sum\limits_{m\neq0}\frac{1}{\langle\sigma+\tfrac{2\pi}{\tau}m\rangle^{2b^\prime}}
\Big|\Big|\sum\limits_{m\neq0}\langle\sigma+\tfrac{2\pi}{\tau}m\rangle^{2b+2b^\prime}\tilde{f}^2(\sigma+\tfrac{2\pi}{\tau}m,k)\Big|\\
&\lesssim\langle\sigma\rangle^{2b}|\tilde{f}(\sigma,k)|^2+\tau^{2b^\prime}\Big|\sum\limits_{m\neq0}
\frac{1}{(\tau\sigma+2\pi m)^{2b^\prime}}\Big|\Big|\sum\limits_{m\neq0}\langle\sigma+\tfrac{2\pi}{\tau}m\rangle^{2b+2b^\prime}\tilde{f}^2(\sigma+\tfrac{2\pi}{\tau}m,k)\Big|\\
&\lesssim\langle\sigma\rangle^{2b}|\tilde{f}(\sigma,k)|^2+\tau^{2b^\prime}\Big|\sum\limits_{m\neq0}\langle\sigma+\tfrac{2\pi}{\tau}m\rangle^{2b+2b^\prime}\tilde{f}^2(\sigma+\tfrac{2\pi}{\tau}m,k)\Big|.
\end{align*}

Integrating the estimate with respect to $\sigma\in[-\tfrac{\pi}{\tau},\tfrac{\pi}{\tau}]$, and taking the square root, we have
$$
\Vert\langle d_\tau\rangle^b\widetilde{f_n}(\cdot,k)\Vert_{L^2\left(-\tfrac{\pi}{\tau},\tfrac{\pi}{\tau}\right)}\lesssim
\Vert\langle\sigma\rangle^b\widetilde{f_n}(\cdot,k)\Vert_{L^2(\mathbb{R})}+\tau^{b^\prime}\Vert\langle\sigma\rangle^{b+b^\prime}\widetilde{f_n}(\cdot,k)\Vert_{L^2(\mathbb{R})}.
$$
The desired inequality then follows from squaring and summing over $k$.
\end{proof}

As a consequence, we obtain the following result.
\begin{proposition}
Let $u^\tau \in X^{s_{0},b_{0}}$ be the extension of the solution of the projected NLS \eqref{p} given by Theorem~\ref{theou-utau} and Remark~\ref{remextension}, and define the sequence $u^\tau_n(x)=u^\tau(n\tau+t^\prime,x)$ for $t^\prime\in[0,4\tau]$.  Then, we have for any $\varepsilon>0$ and $\tau $ sufficiently small the estimate
\begin{equation} \label{h}
\sup\limits_{t^\prime\in[0,4\tau]}\Vert u^\tau_n\Vert_{X_\tau^{s_0,\frac12-\varepsilon}(T)}\leq C_T.
\end{equation}
\end{proposition}

\begin{proof}
By Lemma~\ref{contidisc}, with $b=\frac12-\varepsilon$ for $\varepsilon>0$ and $\frac12 < b'<\frac12+\varepsilon$, we have
$$
\Vert u_n^\tau\Vert_{X_\tau^{s_0,b}} \lesssim\Vert u^{\tau}\Vert_{X^{s_0,b}}+\tau^{b'}\Vert u^{\tau}\Vert_{X^{s_0,b+b'}}.
$$
Consequently, we get from \eqref{bmieux} that
$$
\Vert u_n^\tau\Vert_{X_\tau^{s_0,b}} \lesssim C_{T}+  C_{T}\tau^{b' -{3 \sigma \over 2}}.
$$
Since $b'>\frac12$, we can always choose $\sigma>\frac13$ such that $ b'- {3\sigma\over 2} \geq 0$. This yields
$$
\Vert u_n^\tau\Vert_{X_\tau^{s_0,b}} \lesssim C_{T},
$$
which ends the proof.
\end{proof}

\section{Local error estimates}\label{sectionlocal}

First of all, from the same computations as in  \cite[section~3]{Ost1}, we can
express  the local error
$$
\Psi^\tau(u^\tau(t_n))-u^\tau(t_{n+1})=ie^{i\tau\Delta}\mathcal{E}_{loc}(t_n,\tau,u^\tau),
$$
where $\Psi^\tau$ is defined in \eqref{0} and
\begin{align*}
\mathcal{E}_{loc}(t_n,&\tau,u^\tau)\\
&=\Pi_\tau\int_0^\tau e^{-i\vartheta\Delta}(|\Pi_\tau u^\tau(t_n+\vartheta)|^2\Pi_\tau u^\tau(t_n+\vartheta))+\frac{e^{-i\tau|\Pi_\tau u^\tau(t_n)|^2}-1}{i\tau}\Pi_\tau u^\tau(t_n)d\vartheta\\
&=\Pi_\tau\int_0^\tau\big(e^{-i\vartheta\Delta}-1\big)\big(|\Pi_\tau u^\tau(t_n+\vartheta)|^2\Pi_\tau u^\tau(t_n+\vartheta)\big)d\vartheta\\
&\quad+\Pi_\tau\int_0^\tau|\Pi_\tau u^\tau(t_n+\vartheta)|^2\Pi_\tau\big(u^\tau(t_n+\vartheta)-u^\tau(t_n)\big)d\vartheta\\
&\quad+\Pi_\tau\int_0^\tau\Bigg(|\Pi_\tau u^\tau(t_n+\vartheta)|^2+\frac{e^{-i\tau|\Pi_\tau u^\tau(t_n)|^2}-1}{i\tau}\Bigg)\Pi_\tau u^\tau(t_n)d\vartheta\\
&=\mathcal{E}_1(t_n)+\mathcal{E}_2(t_n)+\mathcal{E}_3(t_n).
\end{align*}
Note that
\begin{equation}\label{g}
\begin{aligned}
u^\tau(t_n+\vartheta)-u^\tau(t_n)&=\big(e^{i\vartheta\Delta}-1\big)u^\tau(t_n)\\
&\quad-i\int_0^\vartheta e^{i(\vartheta-\xi)\Delta}\Pi_\tau\big(|\Pi_\tau u^\tau(t_n+\xi)|^2\Pi_\tau u^\tau(t_n+\xi)\big)d\xi.
\end{aligned}
\end{equation}
We recall that we denote in the following by $u^\tau$ the extension of the solution  of  \eqref{p} on $[0,T]$ defined in Remark~\ref{remextension}. This will not have any influence on the error estimates since we only care about the error on $[0,T]$. In the next theorem, we shall give an estimate of $\mathcal{E}_{loc}(t_n,\tau,u^\tau)$.

\begin{proposition}\label{proplocal}
For $s_0\in(0,2] $ and $u^\tau$ as in Remark~\ref{remextension}, we have for
$\tau$ sufficiently small
$$
\Vert\mathcal{E}_{loc}(t_n,\tau,u^\tau)\Vert_{X^{0,b_0-1}_\tau}\leq C_T\tau^{1+\frac{s_0}{2}},
$$
where $s_0$,  $b_0$ are defined in Theorem~\ref{theoexist}.
\end{proposition}

\begin{proof}
We first estimate $\mathcal{E}_1$. Since $\Pi_\tau$ projects on frequencies less than $\tau^{-\frac{1}{2}}$, we have
$$\sup\limits_{s\in[-\tau,\tau]}\Vert(e^{-is\Delta}-1)\Pi_\tau F(t_n)\Vert_{X^{s,b}_\tau}\lesssim\tau^{\frac{r}{2}}\Vert F(t_n)\Vert_{X^{s+r,b}_\tau}$$
for $r\in[0,2]$, $s,b\in\mathbb{R}$ and any function $F$. Therefore, we get
\begin{align*}
\Vert\mathcal{E}_1(t_n)\Vert_{X^{0,b_{0}-1}_\tau}&\lesssim\left\Vert\Pi_\tau\int_0^\tau(e^{-i\vartheta\Delta}-1)(|\Pi_\tau u^\tau(t_n+\vartheta)|^2\Pi_\tau u^\tau(t_n+\vartheta))\right\Vert_{X^{0,b_{0}-1}_\tau}\\
&\lesssim\tau^{1+\frac{s_{0}}{2}}\sup\limits_{\vartheta\in[0,\tau]}\Vert\Pi_\tau(|\Pi_\tau u^\tau(t_n+\vartheta)|^2\Pi_\tau u^\tau(t_n+\vartheta))\Vert_{X_\tau^{s_{0},b_{0}-1}}\\
&\lesssim\tau^{1+\frac{s_{0}}{2}}\sup\limits_{\vartheta\in[0,\tau]}\Vert u^\tau(t_n+\vartheta)\Vert^3_{X_\tau^{s_{0},1-b_{0}}},
\end{align*}
where the last estimate follows from \eqref{r}. By \eqref{h}, we thus have
\begin{equation}\label{e1}
\Vert\mathcal{E}_1(t_n)\Vert_{X^{0,b_{0}-1}_\tau}\lesssim C_T\tau^{1+\frac{s_{0}}{2}}.
\end{equation}

Next, by \eqref{g}, we get $\mathcal{E}_2(t_n)=\mathcal{E}_{2,1}(t_n)+\mathcal{E}_{2,2}(t_n)$, where
\begin{align*}
&\mathcal{E}_{2,1}(t_n)=\Pi_\tau\int_0^\tau|\Pi_\tau u^\tau(t_n+\vartheta)|^2\Pi_\tau\big((e^{i\vartheta\Delta}-1)u^\tau(t_n)\big)d\vartheta\\
&\mathcal{E}_{2,2}(t_n)=-i\Pi_\tau\int_0^\tau|\Pi_\tau u^\tau(t_n+\vartheta)|^2\Pi_\tau\int_0^\vartheta e^{i(\vartheta-\xi)\Delta}\Pi_\tau\big(|\Pi_\tau u^\tau(t_n+\xi)|^2\Pi_\tau u^\tau(t_n+\xi)\big)d\xi d\vartheta.
\end{align*}
First, by using \eqref{rbis}, we have
\begin{align*}
\Vert\mathcal{E}_{2,1}(t_n)\Vert_{X^{0,b_{0}-1}_\tau}&\lesssim\tau\sup\limits_{\vartheta\in[0,\tau]}\Big(\left\Vert\Pi_\tau\Bigl(|\Pi_\tau u^\tau(t_n+\vartheta)|^2\Pi_\tau\big((e^{i\vartheta\Delta}-1)u^\tau(t_n)\big)\Bigr)\right\Vert_{X_\tau^{0,b_{0}-1}}\Big)\\
&\lesssim\tau\sup\limits_{\vartheta\in[0,\tau]}\Big(\Vert u^\tau(t_n+\vartheta)\Vert^2_{X^{s_{0},1-b_{0}}_\tau}\Vert(e^{i\vartheta\Delta}-1)\Pi_\tau u^\tau(t_n)\Vert_{X^{0,1-b_{0}}_\tau}\Big)\\
&\lesssim\tau^{1+\frac{s_{0}}{2}}\sup\limits_{\vartheta\in[0,\tau]}\Vert u^\tau(t_n+\vartheta)\Vert^3_{X_\tau^{s_{0},1-b_{0}}}.
\end{align*}
Again, by \eqref{h}, we thus obtain
\begin{equation}\label{e2}
\Vert\mathcal{E}_{2,1}(t_n)\Vert_{X^{0,b_{0}-1}_\tau}\lesssim C_T\tau^{1+\frac{s_{0}}{2}}.
\end{equation}
By using again \eqref{rbis} and \eqref{a},  we get that
\begin{equation} \label{E22debut}
\begin{aligned}
\Vert\mathcal{E}_{2,2}&(t_n)\Vert_{X^{0,b_{0}-1}_\tau}\\
&\lesssim\tau^2\sup\limits_{\vartheta\in[0,\tau]}\Big(\Vert \Pi_{\tau }u^\tau(t_n+\vartheta)\Vert_{X^{s_{0},1-b_{0}}_\tau}^2 \Vert\Pi_{\tau} \left( |\Pi_\tau  u^\tau(t_n+\vartheta)|^2\Pi_\tau u^\tau(t_n+\vartheta) \right)\Vert_{X_\tau^{0,1-b_{0}}}\Big).
\end{aligned}
\end{equation}
From \eqref{b}, we then  get that
\begin{equation}\label{E22debut2}
\tau\Vert\Pi_{\tau} \left( |\Pi_\tau  u^\tau(t_n+\vartheta)|^2\Pi_\tau u^\tau(t_n+\vartheta) \right)\Vert_{X_\tau^{0,1-b_{0}}} \lesssim \tau^{b_{0}}\Vert\Pi_{\tau} \left( |\Pi_\tau  u^\tau(t_n+\vartheta)|^2\Pi_\tau u^\tau(t_n+\vartheta) \right)\Vert_{X_\tau^{0,0}}.
\end{equation}
Next, since $X_\tau^{0,0}=l^2_{\tau}L^2, $ we actually have
\begin{equation}\label{E22+}
\tau^{b_{0}}\Vert\Pi_{\tau} \left( |\Pi_\tau  u^\tau(t_n+\vartheta)|^2\Pi_\tau u^\tau(t_n+\vartheta) \right)\Vert_{X_\tau^{0,0}} \lesssim   \tau^{b_{0}} \| \Pi_{\tau} u^\tau(t_{n}+ \vartheta) \|_{l^6_{\tau}L^6}^3.
\end{equation}
By using the Sobolev embedding $W^{{1\over 6}, 4}(\mathbb{T}^2) \subset L^6 $ and that  for every sequence $(f_{n})_{n}$, we have that
\begin{equation}\label{l6l4}
\|f_{n}\|_{l^6_{\tau}} \lesssim \tau^{-\frac1{12}} \|f_{n}\|_{l^4_{\tau}},
\end{equation}
we then obtain
$$
\tau^{b_{0}} \| \Pi_{\tau} u^\tau(t_{n}+ \vartheta) \|_{l^6_{\tau}L^6}^3
\lesssim   \tau^{b_{0}- {1 \over 4}} \| \langle \partial_{x} \rangle^{1\over 6}\Pi_{\tau} u^\tau (t_{n}+ \vartheta) \|_{l^4_{\tau}L^4}^3.
$$
Consequently, by using the Strichartz estimate \eqref{z}, we get
\begin{equation}\label{E22milieu}
\tau^{b_{0}} \| \Pi_{\tau} u^\tau(t_{n}+ \vartheta) \|_{l^6_{\tau}L^6}^3\ \lesssim   \tau^{b_{0}- {1 \over 4}}
\| \Pi_{\tau} u^\tau(t_{n}+ \vartheta) \|_{X_\tau^{{1 \over 6} + \varepsilon, {1 \over 2} + \varepsilon}}^3,
\end{equation}
where $\varepsilon>0$ can be taken arbitrarily small.

Thus, if $s_{0} \leq {1 \over 6}$, we get by using again \eqref{b}, \eqref{t} that
\begin{equation}\label{presqueE3}
\begin{aligned}
\tau^{b_{0}} \| \Pi_{\tau} u^\tau(t_{n}+ \vartheta) \|_{l^6_{\tau}L^6}^3 &\lesssim   \tau^{b_{0}- {1 \over 4} - {3 \over 2}( { 1\over 6}+\varepsilon-s_{0})- 3({1 \over 2} + \varepsilon- (1-b_{0}))}
\| \Pi_{\tau} u^\tau(t_{n}+ \vartheta) \|_{X_\tau^{s_{0}, 1-b_{0}}}^3\\
&\lesssim  \tau^{ -2b_{0} + 1 + {3 \over 2} s_{0} - { 9 \over 2} \varepsilon} \|  \Pi_{\tau} u^\tau(t_{n}+ \vartheta) \|_{X_\tau^{s_{0}, 1-b_{0}}}^3.
\end{aligned}
\end{equation}
This yields thanks to \eqref{E22debut} and \eqref{h} again that
$$
\Vert\mathcal{E}_{2,2}(t_n)\Vert_{X^{0,b_{0}-1}_\tau} \lesssim C_{T} \tau^{1 -2b_{0} + 1 + {3 \over 2} s_{0} - { 9 \over 2} \varepsilon} \lesssim C_{T} \tau^{ 1 + {s_{0}\over 2}}.
$$
Indeed, since $\varepsilon$ can be taken arbitrarily small, we just need
$$
-2b_{0} + 1 + {3 \over 2} s_{0}>{s_{0}\over 2},
$$
which is equivalent to
$$
b_{0}<{1 \over 2} + {s_{0} \over 2}.
$$
This is always satisfied when $b_{0}$ is taken as in Theorem~\ref{theoexist}.

It remains the case $s_{0}>{1 \over 6}$. In this case, we get from \eqref{E22milieu} that
\begin{equation}\label{presqueE32}
\begin{aligned}
\tau^{b_{0}} \| \Pi_{\tau} u^\tau(t_{n}+ \vartheta) \|_{l^6_{\tau}L^6}^3&\lesssim  \tau^{b_{0}- {1 \over 4} - 3({1 \over 2} + \varepsilon- (\frac12-\varepsilon))}
 \| \Pi_{\tau} u^\tau (t_{n}+ \vartheta) \|_{X_\tau^{s_{0}, \frac12-\varepsilon}}^3\\
&\lesssim \tau^{b_{0}- {1\over 4} - 6 \varepsilon} \|  \Pi_{\tau} u^\tau (t_{n}+ \vartheta) \|_{X_\tau^{s_{0}, \frac12-\varepsilon}}^3,
\end{aligned}
\end{equation}
which yields, thanks to \eqref{E22debut} and \eqref{h},  that
$$
\Vert\mathcal{E}_{2,2}(t_n)\Vert_{X^{0,b_{0}-1}_\tau} \lesssim C_{T} \tau^{1+b_{0}- {1\over 4} - 6\varepsilon} \lesssim  C_{T}  \tau^{1 + {s_{0} \over 2}}
$$
since it suffices to ensure that
$$
b_{0} - {1\over 4}>{s_{0}\over 2}.
$$
Since $b_{0}>\frac12$, this condition is satisfied for $s_0<\frac12$.

It thus remains to consider $s_0 \geq \frac12$. For this case, we start with the crude estimate
\begin{equation}\label{presqueE33}
\Vert\mathcal{E}_{2,2}(t_n)\Vert_{X^{0,b_{0}-1}_\tau}\lesssim \Vert\mathcal{E}_{2,2}(t_n)\Vert_{l^2_{\tau}L^2} \lesssim \tau^2 \sup_{\vartheta \in [0, \tau]} \| u^\tau(t_n+\vartheta) \|_{l^{10}_{\tau}L^{10}}^5.
\end{equation}
We now use the Sobolev embedding $H^{4\over 5}(\mathbb{T}^2) \subset L^{10}$ and the fact that the sequence $u(t_{n}+ \vartheta)$ is compactly supported in $n \tau \lesssim 1$. This yields
\begin{equation} \label{presquefini}
\Vert\mathcal{E}_{2,2}(t_n)\Vert_{X^{0,b_{0}-1}_\tau}\lesssim \Vert\mathcal{E}_{2,2}(t_n)\Vert_{l^2_{\tau}L^2} \lesssim \tau^2  \sup_{\vartheta \in [0, \tau]} \|\Pi_{\tau} u^\tau(t_n+\vartheta) \|_{l^{\infty}_{\tau}H^{4\over5}}^5 \lesssim \tau^2 \|\Pi_{\tau}u^\tau \|_{L^\infty H^{4\over 5}}^5.
\end{equation}
We then have again two cases. If $\frac12 \leq s_{0} \leq \frac45$, thanks to the frequency localization induced by $\Pi_{\tau}$, we write
$$
\Vert\mathcal{E}_{2,2}(t_n)\Vert_{X^{0,b_{0}-1}_\tau}\lesssim  \tau^{2 - {5 \over 2}({4\over 5}- s_{0})}\|u^\tau \|_{L^\infty H^{s_{0}}}^5\lesssim \tau^{{ 5 \over 2}s_{0}} C_{T},
$$
where we have used the embedding $X^{s_{0}, b_{0}} \subset L^\infty H^{s_{0}}$ for the last step. If ${ 5 \over 2} s_{0} \geq 1+ {s_{0} \over 2}$, which is equivalent to $s_{0} \geq {1\over 2}$,
we thus get
$$
\Vert\mathcal{E}_{2,2}(t_n)\Vert_{X^{0,b_{0}-1}_\tau}\lesssim \tau^{1 + {s_{0}\over 2}} C_{T}.
$$
If $s_{0}>\frac45$, then we directly get from \eqref{presquefini} and since $X^{s_{0}, b_{0}} \subset L^\infty H^{s_{0}}$ that
$$
\Vert\mathcal{E}_{2,2}(t_n)\Vert_{X^{0,b_{0}-1}_\tau}\lesssim  C_{T} \tau^2 \lesssim   C_{T} \tau^{ 1 + {s_{0}\over 2}}
$$
since  $s_{0} \leq 2$.
In summary, we have finally obtained that for all $s_{0}\in (0, 2]$, we have
$$
\Vert\mathcal{E}_{2,2}(t_n)\Vert_{X^{0,b_{0}-1}_\tau}\lesssim C_{T}\tau^{ 1 + {s_{0}\over 2}}.
$$
Gathering the estimates for $\mathcal{E}_{2,2}$ and \eqref{e2}, we finally obtain the estimate of $\mathcal{E}_{2}$,
\begin{equation}\label{e2fin}
\Vert\mathcal{E}_{2}(t_n)\Vert_{X^{0,b_{0}-1}_\tau}\lesssim \tau^{1 + {s_{0}\over 2}} C_{T}
\end{equation}
for all $s_{0} \in (0, 2]$.

To estimate $\mathcal{E}_3$, we rewrite it as
$$
\mathcal{E}_3(t_n)=\mathcal{E}_{3,1}(t_n)+\mathcal{E}_{3,2}(t_n)
$$
with
\begin{align*}
\mathcal{E}_{3,1}(t_n)&=\Pi_\tau\int_0^\tau\big(|\Pi_\tau u^\tau(t_n+\vartheta)|^2-|\Pi_\tau u^\tau(t_n)|^2\big)\Pi_\tau u^\tau(t_n)d\vartheta\\
&= \Pi_\tau\int_0^\tau\big(\Pi_\tau(u^\tau(t_n+\vartheta)-u^\tau(t_n))\big) \cdot \big(\Pi_\tau u^\tau(t_n+\vartheta)+\Pi_\tau u^\tau(t_n)\big)\Pi_\tau u^\tau(t_n)d\vartheta,\\
\mathcal{E}_{3,2}(t_n)&=\tau\Pi_\tau\Bigg(\dfrac{e^{-i\tau|\Pi_\tau u^\tau(t_n)|^2}-1+i\tau|\Pi_\tau u^\tau(t_n)|^2}{i\tau}\Pi_\tau u^\tau(t_n)\Bigg).
\end{align*}

By using the same arguments as for the estimate of $\mathcal{E}_2(t_n)$, we find that
\begin{equation}\label{e4}
\Vert\mathcal{E}_{3,1}(t_n)\Vert_{X^{0,b_0-1}_\tau}\lesssim C_T\tau^{1+\frac{s_{0}}{2}}.
\end{equation}
The estimate of  $\mathcal{E}_{3,2}$ will be also rather similar to the one of  $\mathcal{E}_{2,2}$ since
\begin{equation}\label{expfacile}
\Big|\dfrac{e^{-i\tau\alpha}-1+i\tau\alpha}{i\tau}\Big|\lesssim\tau|\alpha|^2,\quad \forall \alpha\in\mathbb{R}.
\end{equation}
We first rewrite $\mathcal{E}_{3, 2}$ in the form
$$
\mathcal{E}_{3,2}(t_n) =\tau\Pi_\tau\Bigg( | \Pi_{\tau}u^{\tau}(t_{n})|^2 \left( \, \dfrac{e^{-i\tau|\Pi_\tau u^\tau(t_n)|^2}-1+i\tau|\Pi_\tau u^\tau(t_n)|^2}{i\tau |\Pi_{\tau }u^\tau (t_{n})|^2 } \Pi_\tau u^\tau(t_n) \right)\Bigg).
$$
Consequently, by using \eqref{rbis} again, we get
$$
\Vert\mathcal{E}_{3,2}(t_n)\Vert_{X^{0,b_0-1}_\tau}\lesssim\tau\|\Pi_{\tau}u^\tau(t_n)\|_{X^{s_{0}, 1-{b_{0}}}_{\tau}}^2
\left\| \dfrac{e^{-i\tau|\Pi_\tau u^\tau(t_n)|^2}-1+i\tau|\Pi_\tau u^\tau(t_n)|^2}{i\tau |\Pi_{\tau }u^\tau (t_{n})|^2 } \Pi_\tau u^\tau(t_n)\right\|_{X^{0, 1-b_{0}}_{\tau}}.
$$
To estimate the right-hand side, we use \eqref{b} again to get
\begin{multline*}
\left\| \dfrac{e^{-i\tau|\Pi_\tau u^\tau(t_n)|^2}-1+i\tau|\Pi_\tau u^\tau(t_n)|^2}{i\tau |\Pi_{\tau }u^\tau (t_{n})|^2 } \Pi_\tau u^\tau(t_n)\right\|_{X^{0, 1-b_{0}}_{\tau}} \\ \lesssim \tau^{b_{0}-1}
\left\| \dfrac{e^{-i\tau|\Pi_\tau u^\tau(t_n)|^2}-1+i\tau|\Pi_\tau u^\tau(t_n)|^2}{i\tau |\Pi_{\tau }u^\tau (t_{n})|^2 } \Pi_\tau u^\tau(t_n)\right\|_{X^{0,0}_{\tau}}
\end{multline*}
and hence, thanks to \eqref{expfacile}, we have
$$
\left\| \dfrac{e^{-i\tau|\Pi_\tau u^\tau(t_n)|^2}-1+i\tau|\Pi_\tau u^\tau(t_n)|^2}{i\tau |\Pi_{\tau }u^\tau (t_{n})|^2 } \Pi_\tau u^\tau(t_n)\right\|_{X^{0, 1-b_{0}}_{\tau}} \lesssim \tau^{b_{0}}
\|\Pi_{\tau}u^\tau(t_{n})\|_{l^6_{\tau}L^6}^3.
$$
Consequently, by using \eqref{presqueE3}, \eqref{presqueE32}, we also get that
$$
\Vert\mathcal{E}_{3,2}(t_n)\Vert_{X^{0,b_0-1}_\tau}\lesssim  C_{T}\tau^{ 1 + { s_{0} \over 2}}
$$
for $s_{0} <\frac12$.

For the case $s_{0} \geq \frac12$, we just use  that
$$
\Vert\mathcal{E}_{3,2}(t_n)\Vert_{X^{0,b_0-1}_\tau}\lesssim \Vert\mathcal{E}_{3,2}(t_n)\Vert_{X^{0,0}_\tau}\lesssim \tau^2 \|u^\tau(t_{n})\|_{l^{10}_{\tau}L^{10}}^5
$$
thanks to \eqref{expfacile}. We thus have the same upper bound as in \eqref{presqueE33} and we can use the above estimates to get that
$$
\Vert\mathcal{E}_{3,2}(t_n)\Vert_{X^{0,b_0-1}_\tau}\lesssim   C_{T}\tau^{ 1 + { s_{0} \over 2}}
$$
for ${1 \over 2} \leq s_{0} \leq 2$, so that the estimate
\begin{equation} \label{e3}
\Vert\mathcal{E}_{3,2}(t_n)\Vert_{X^{0,b_0-1}_\tau}\lesssim   C_{T}\tau^{ 1 + { s_{0} \over 2}}
\end{equation}
holds for all $s_{0}\in (0, 2].$

We finish the proof by collecting \eqref{e1}, \eqref{e2fin}, \eqref{e4}, \eqref{e3}.
\end{proof}

\section{Global error estimates}\label{sectionglobal}

First of all, similar to \cite[section~3]{Ost1}, we can write the global error as
\begin{equation}\label{k}
\begin{aligned}
e_n&=u^\tau(t_n)-u_n\\
&=u^\tau(t_{n-1})-u_{n-1}-i\tau e^{i\tau\Delta}(\Phi^\tau(u^\tau(t_{n-1}))-\Phi^\tau(u_{n-1}))-ie^{i\tau\Delta}\mathcal{E}_{loc}(t_{n-1},\tau,u^\tau)\\
&=-i\tau\sum\limits_{k=0}^{n-1}e^{i(n-k)\tau\Delta}(\Phi^\tau(u^\tau(t_k))-\Phi^\tau(u_k))-i\sum\limits_{k=0}^{n-1}e^{i(n-k)\tau\Delta}\mathcal{E}_{loc}(t_k,\tau,u^\tau)
\end{aligned}
\end{equation}
with the nonlinear flow
\begin{equation}\label{defphitaufin}
\Phi^\tau(w)=-\Pi_\tau\Big(\dfrac{e^{-i\tau|\Pi_\tau w|^2}-1}{i\tau}\Pi_\tau w\Big).
\end{equation}

In this section, we estimate the global error $e_n$.
\begin{proposition}\label{propen}
For $s_{0} \in (0,2]$ and $\tau$ sufficiently small, we have
$$
\Vert e_n\Vert_{X^{0,b_{0}}_{\tau}(T)} \lesssim \tau^{ s_{0}\over 2},
$$
where $s_0$, $b_0$ are defined in Theorem~\ref{theoexist}.
\end{proposition}

\begin{proof}
Take a smooth function $\eta$ which is 1 on $[-1,1]$ and compactly supported in $[-2,2]$. In the proof we shall still denote the solution of the truncated version of \eqref{k} by $e_{n}$, given as
\begin{equation}\label{m}
e_n=-i\tau\eta(t_n)\sum\limits_{k=0}^{n-1}e^{i(n-k)\tau\Delta}\,\eta\biggl(\frac{t_k}{T_1}\biggr)\big(\Phi^\tau(u^\tau(t_k))-\Phi^\tau(u^\tau(t_k)-e_{k})\big)+\mathcal{R}_n,
\end{equation}
where
$$
\mathcal{R}_n=-i\eta(t_n)\sum\limits_{k=0}^{n-1}e^{i(n-k)\tau\Delta}\eta\biggl(\frac{t_k}{T_1}\biggr)\mathcal{E}_{loc}(t_k,\tau,u^\tau).
$$
Note that for $0\leq n\leq N_1$, where $N_1=[\frac{T_1}{\tau}]$ with $T_1\leq\min(1,T)$, this indeed coincides with $u^\tau(t_{n})- u_{n}$. By using property \eqref{q} which was stated in Lemma~\ref{discprop}, we get
$$
\Vert\mathcal{R}_n\Vert_{X^{0,b_{0}}_\tau}\lesssim\tau^{-1}\Vert\mathcal{E}_{loc}(t_n,\tau,u^\tau)\Vert_{X^{0,b_{0}-1}_\tau}.
$$
By Proposition~\ref{proplocal}, we immediately get
\begin{equation}\label{8}
\Vert\mathcal{R}_n\Vert_{X^{0,b_{0}}_\tau}\leq C_T\tau^{\frac{s_{0}}{2}}.
\end{equation}
Thus by  \eqref{m} and \eqref{8}, we have
$$
\Vert e_n\Vert_{X_{\tau}^{0,b_{0}}}\lesssim\left\Vert\tau\eta(t_n)\sum\limits_{k=0}^{n-1}e^{i(n-k)\tau\Delta}\eta\left(\frac{t_k}{T_1}\right)\big(\Phi^\tau(u^\tau(t_k))-\Phi^\tau(u^\tau(t_k)-e_k)\big)\right\Vert_{X_{\tau}^{0,b_{0}}}+
C_T\tau^{\frac{s_{0}}{2}}.
$$
By the definition of $\Phi^\tau$, see \eqref{defphitaufin}, we can write it in the form
$$
\Phi^{\tau}= \Pi_{\tau}F^{\tau}(w) + \Pi_{\tau} R^{\tau}(w),
$$
where
\begin{align*}
F^{\tau}(w)&= |\Pi_{\tau}(w)|^2 \Pi_{\tau}w, \\
R^{\tau}(w)& = F^\tau(w) R^{\tau, 1}(w), \quad R^{\tau, 1}(w)= \int_{0}^1  \left( e^{-is \tau  |\Pi_{\tau}w|^2}- 1\right) ds.
\end{align*}
Note that $R^{\tau}(w)$ behaves like a quintic term in the sense that
$$
|R^{\tau}(w)| \lesssim \tau |\Pi_{\tau}w |^5
$$
with the gain of a factor $\tau$. The main idea of this decomposition is the following: we use the multilinear estimate \eqref{rbis} to estimate the leading order cubic term exactly, since we cannot lose anything in its estimation. For the other term, because of the gain of the factor $\tau$, we can use the Strichartz estimate \eqref{z} instead to estimate it, thus taking advantage of the fact that the exponential factor has modulus one.

By using this decomposition and \eqref{9}, \eqref{q}, we get that
\begin{equation} \label{esten1}
\begin{aligned}
\Vert e_n\Vert_{X^{0,b_{0}}_{\tau}}\leq C_TT_1^{\varepsilon_0}\Vert \Pi_{\tau}&\bigl(F^\tau(u^\tau(t_n))-F^\tau(u^\tau(t_n)-e_n) \bigr)\Vert_{X_\tau^{0,-b'}} \\
&+ C_{T}\Vert\Pi_{\tau}\bigl(R^{\tau}(u^\tau(t_n))-R^{\tau}(u^\tau(t_n)-e_n) \bigr)\Vert_{X_\tau^{0,b_{0}-1}} +C_T\tau^{\frac{s_{0}}{2}},
\end{aligned}
\end{equation}
where $b'$ is chosen as in the proofs of Theorem~\ref{theoexist} and Theorem~\ref{theou-utau}. That is to say, it satisfies $\max(\frac{1}{4},\frac{1}{2}-\frac{1}{4}s_0)<b'<1 - b_{0}$, which implies that $\varepsilon_0= 1- b_{0}-b'>0$.

Since we can expand
$$
\Pi_{\tau}\Bigl(F^\tau(u^\tau(t_n))-F^\tau(u^\tau(t_n)-e_n) \Bigr)= L^\tau_n + Q^\tau_n  + C^\tau_n
$$
with
\begin{align*}
L^\tau_{n}&=\Pi_{\tau} \left( |\Pi_{\tau} u^\tau(t_{n})|^2\Pi_{\tau} e_{n} + 2 \Pi_{\tau} u^\tau(t_{n}) \cdot  \Pi_{\tau}e_{n} \Pi_{\tau} u^\tau(t_{n}) \right), \\
Q^{\tau}_{n}&=  -  \Pi_{\tau} \left( 2 \Pi_{\tau} u^{\tau}(t_{n}) \cdot \Pi_{\tau} e_{n} \Pi_{\tau}e_{n} + \Pi_{\tau}u^\tau(t_{n})|\Pi_{\tau} e_{n}|^2\right), \\
C^\tau_{n}&= \Pi_{\tau}\left( |\Pi_{\tau} e_{n}|^2 \Pi_{\tau}e_{n} \right),
\end{align*}
we get by using \eqref{rbis} that
\begin{equation}\label{ouf4}
\begin{aligned}
\Vert \Pi_{\tau}\bigl(F^\tau(u^\tau(t_n))&-F^\tau(u^\tau(t_n)-e_n) \bigr)\Vert_{X_\tau^{0,-b'}}
\lesssim \|u^\tau(t_{n}) \|_{X^{s_{0}, 1-b_{0}}_{\tau}}^2 \| e_{n}\|_{X^{0, 1 -b_{0}}_{\tau}} \\
&+  \|u^\tau(t_{n}) \|_{X^{s_{0}, 1-b_{0}}_{\tau}} \| e_{n}\|_{X^{0, 1 -b_{0}}_{\tau}} \|\Pi_{\tau}e_{n}\|_{ X^{s_1, 1 -b_{0}}_{\tau}}
+\|\Pi_{\tau}e_{n}\|^2 _{X^{s_1, 1 -b_{0}}_{\tau}} \| e_{n}\|_{X^{0, 1 -b_{0}}_{\tau}},
\end{aligned}
\end{equation}
with $s_1\in (0, s_{0}]$ to be chosen. Note that we have freely added $\Pi_{\tau}$ in the second term and the third term of the right-hand side, which is allowed since $\Pi_{\tau} e_{n} = e_{n}$. By using \eqref{t}, which reads
$$
\|\Pi_{\tau}e_{n}\|_{X^{s_1, 1 -b_{0}}_{\tau}} \lesssim \tau^ {- {s_1 \over 2}} \|e_{n}\|_{X^{0, 1 -b_{0}}_{\tau}},
$$
and \eqref{h}, we thus get that
\begin{equation}\label{esten2}
\begin{aligned}
\Vert \Pi_{\tau}\bigl(F^\tau(u^\tau(t_n))-F^\tau(&u^\tau(t_n)-e_n) \bigr)\Vert_{X_\tau^{0,-b'}} \\
& \leq C_{T}  \|e_{n}\|_{X^{0, 1 -b_{0}}_{\tau}} + C_{T,s_1 } \tau^{- {s_1 \over 2}}  \|e_{n}\|_{X^{0, 1 -b_{0}}_{\tau}}^2
+ C_{T, s_1} \tau^{-  s_1}  \|e_{n}\|_{X^{0, 1 -b_{0}}_{\tau}}^3.
\end{aligned}
\end{equation}
We shall now estimate $\Vert \Pi_{\tau}\left(R^{\tau}(u^\tau(t_n))-R^{\tau}(u^\tau(t_n)-e_n) \right)\Vert_{X_\tau^{0,b_0-1}}$. We first use \eqref{b} and \eqref{t} to get
\begin{multline*}
\Vert \Pi_{\tau}\bigl(R^{\tau}(u^\tau(t_n))-R^{\tau}(u^\tau(t_n)-e_n) \bigr)\Vert_{X_\tau^{0,b_{0}-1}}\\
\lesssim \tau^{-\frac{3\delta}2 - (b_{0} -\frac12)}  \Vert \Pi_{\tau}\bigl(R^{\tau}(u^\tau(t_n))-R^{\tau}(u^\tau(t_n)-e_n) \bigr)\Vert_{X_\tau^{-\delta,-{1 \over 2}- \delta}},
\end{multline*}
where $\delta >0$ can be arbitrarily small and will be chosen small enough. Now we can use the dual version of the Strichartz estimate \eqref{z}, which reads for any $(u_{n})$,
$$
\|\Pi_{\tau} u_{n} \|_{X^{-\delta, -{1 \over 2}- \delta}} \lesssim  \|u_{n}\|_{l^{4\over 3}_{\tau}L^{4\over3}}.
$$
This yields
$$
\Vert \Pi_{\tau}\left(R^{\tau}(u^\tau(t_n))-R^{\tau}(u^\tau(t_n)-e_n) \right)\Vert_{X_\tau^{0,b_0-1}}
\lesssim \tau^{-\frac{3\delta}2 - (b_{0} -\frac12)}  \Vert R^{\tau}(u^\tau(t_n))-R^{\tau}(u^\tau(t_n)-e_n) \Vert_{l^{4\over3}_{\tau}L^{4\over3}}.
$$
Now, we observe that we have the pointwise uniform in $\tau$ estimate
$$
\left| \left(R^{\tau}(u^\tau(t_n))-R^{\tau}(u^\tau(t_n)-e_n) \right) \right|  \lesssim \tau\sum_{j=1}^5 |u^\tau(t_{n})|^{5-j} |e_{n}|^j .
$$
We thus deduce that
\begin{equation}\label{ouf3}
\Vert \Pi_{\tau}\left(R^{\tau}(u^\tau(t_n))-R^{\tau}(u^\tau(t_n)-e_n) \right)\Vert_{X_\tau^{0,b_0-1}}
\lesssim \tau^{1-\frac{3\delta}2 -( b_{0} -\frac12)}
\sum_{j=1}^5 \left\| |\Pi_{\tau}u^\tau(t_{n})|^{5-j} |\Pi_{\tau}e_{n}|^j \right\|_{l^{4\over 3}_{\tau}L^{4\over3}}.
\end{equation}
From H\"{o}lder's inequality, we get that
\begin{align*}
\sum_{j=1}^5 &\left\| |\Pi_{\tau }u^\tau(t_{n})|^{5-j} |\Pi_{\tau}e_{n}|^j \right\|_{l^{4\over 3}_{\tau}L^{4\over3}} \lesssim
\|\Pi_{\tau}u^\tau(t_{n})\|_{l^\infty_{\tau}L^\infty}^2 \|\Pi_{\tau}u^\tau(t_{n})\|_{l^4_{\tau}L^4}^2 \|\Pi_{\tau} e_{n}\|_{l^4_{\tau}L^4}\\[-2.5mm]
&\qquad\qquad + \|\Pi_{\tau}u^\tau(t_{n})\|_{l^\infty_{\tau}L^\infty}^2 \|\Pi_{\tau}u^\tau(t_{n})\|_{l^4_{\tau}L^4} \|\Pi_{\tau} e_{n}\|_{l^4_{\tau}L^4}^2
+ \|\Pi_{\tau}u^\tau(t_{n})\|_{l^\infty_{\tau}L^\infty}^2  \|\Pi_{\tau} e_{n}\|_{l^4_{\tau}L^4}^3\\
&\qquad\qquad +
\|\Pi_{\tau}u^\tau(t_{n})\|_{l^\infty_{\tau}L^\infty} \|\Pi_{\tau}e_{n}\|_{l^\infty_{\tau}L^\infty}  \|\Pi_{\tau} e_{n}\|_{l^4_{\tau}L^4}^3 +   \|\Pi_{\tau}e_{n}\|_{l^\infty_{\tau}L^\infty}^2 \|\Pi_{\tau} e_{n}\|_{l^4_{\tau}L^4}^3.
\end{align*}
To estimate the right-hand side, we use the  Strichartz estimate \eqref{z}, and again \eqref{h} and \eqref{t} to get
$$
\|\Pi_{\tau} e_{n}\|_{l^4_{\tau}L^4} \lesssim \|\Pi_{\tau}e_{n}\|_{X^{\delta, {1 \over 2}+\delta}_{\tau}}
\lesssim \tau^{-\frac{\delta}2} \|e_{n}\|_{X^{0, b_{0}}_{\tau}}
$$
by choosing $\delta$ sufficiently small so that $b_{0}>1/2 + \delta$. In a similar way, we get that
$$
\|\Pi_\tau u^\tau(t_n)\|_{l^4_{\tau}L^4} \lesssim \|\Pi_\tau u^\tau(t_n)\|_{X^{\delta,{1\over2}+\delta}_{\tau}}
\lesssim \tau^{-2{\delta}} \|u^\tau(t_{n})\|_{X^{s_{0},{1 \over 2} - \delta}_{\tau}}
\lesssim  \tau^{-2{\delta}}  C_{T}
$$
by choosing $\delta <s_{0}$. The last estimate comes from \eqref{h}.

Moreover, by Sobolev embedding, the continuity of $u^\tau$ and \eqref{d}, we can also get that
\begin{align*}
\| \Pi_{\tau}u^\tau(t_{n})\|_{l^\infty_{\tau}L^\infty}&
\lesssim  \tau^{-{1 \over 2}(1+ \delta - s_0)}\|u^\tau(t_n)\|_{l^\infty_{\tau}H^{s_{0}}} \lesssim \tau^{-{1 \over 2}(1+\delta - s_0)} \|u^\tau\|_{L^\infty H^{s_0}}\\
&\lesssim \tau^{-{1 \over 2}(1+\delta - s_0)}\|u^\tau\|_{X^{s_{0},b_0}}\leq\tau^{-{1 \over 2}(1+\delta - s_0)} C_T
\end{align*}
if $s_{0} \leq 1$. In the case $s_{0}>1$, we just have a simpler estimate directly from the Sobolev embedding and a similar argument
$$
\| \Pi_{\tau}u^\tau(t_{n})\|_{l^\infty_{\tau}L^\infty} \lesssim  {\|u^\tau\|_{L^\infty H^{s_{0}}}\lesssim C_{T}}.
$$
Finally, we also have that
$$
\| \Pi_{\tau} e_{n}\|_{l^\infty_{\tau}L^\infty}
\lesssim  \tau^{-{1 \over 2} - {\delta \over 2}} \|e_{n}\|_{l^\infty_{\tau}L^2}
\lesssim   \tau^{-{1 \over 2} - {\delta \over 2}} \|e_{n}\|_{X^{0, b_{0}}_{\tau}},
$$
where we have used again \eqref{y} for the last part of the estimate.

This yields in the case $s_{0} \leq 1$ (the case $s_{0}>1$ is much easier and can be handled by similar arguments),
\begin{equation}\label{ouf}
\begin{aligned}
\tau \sum_{j=1}^5 \left\| |\Pi_{\tau }u^\tau(t_{n})|^{5-j} |\Pi_{\tau}e_{n}|^j \right\|_{l^{4\over 3}_{\tau}L^{4\over3}}
&\lesssim C_{T,\delta} \left(   \tau^{ s_{0} -  {11\delta \over 2}}
\|e_{n}\|_{X_{\tau}^{0, b_{0}}} +  \tau^{ s_{0} - 4 \delta} \|e_{n}\|_{X_{\tau}^{0, b_{0}}}^2 \right. \\
&\hspace{-5mm}\left. {}+ { \tau^{ s_{0}-  { 5\delta \over 2}}}  \|e_{n}\|_{X_{\tau}^{0, b_{0}}}^3+ \tau^{ {s_{0}\over 2}- {5 \delta \over 2}} \|e_{n}\|_{X_{\tau}^{0, b_{0}}}^4 +\tau^{-{5 \delta \over 2} }  \|e_{n}\|_{X_{\tau}^{0, b_{0}}}^5 \right).
\end{aligned}
\end{equation}
Consequently, by combining \eqref{esten1}, \eqref{ouf4}, \eqref{ouf3}, and \eqref{ouf}, we obtain that
\begin{multline*}
\|e_{n}\|_{X^{0, b_0}_\tau} \leq C_{T} \left( ( T_1^{\varepsilon_0} + C_{T,\delta}
\tau^{s_0 -7\delta - (b_0- {1\over 2})}) \|e_n\|_{X^{0, b_0}_\tau} \right. \\
\left.{} + C_{T, s_1,\delta} (\tau^{-{s_1\over 2}} + \tau^{s_0- {11\delta \over 2} - (b_0-{1\over2})})
\|e_n\|_{X^{0, b_0}_\tau}^2 + C_{T, s_1,\delta} ( \tau^{-s_1} +\tau^{s_{0} -  4 \delta  - (b_{0}- {1 \over 2})}) \|e_n\|_{X^{0, b_0}_\tau}^3 \right. \\
\left.{} + C_{T, \delta}\tau^{{s_0\over 2}- 4 \delta - (b_0- {1\over2})}\|e_n\|_{X^{0, b_0}_\tau}^4
+  C_{T, \delta} \tau^{-4\delta - (b_0-{1\over 2})} \|e_n\|_{X^{0,b_0}_\tau}^5 + \tau^{s_0\over 2}\right).
\end{multline*}
To conclude, we observe that for $s_0>0$ fixed, we can choose $\delta$ sufficiently small ($s_0> 50\delta$ for example) and $s_1$ satisfying $1-b_0>b^\prime>\tfrac12-\tfrac14s_1$ and $\tfrac{s_0}2-\tfrac{s_1}2>0$, i.e. $4b_0-2<2-4b^\prime<s_1<s_0$. This is possible since we always have $b^\prime>\tfrac12-\tfrac14s_0$. Since $b_0-\tfrac12<\tfrac14s_0$, we can get
\begin{multline*}
\|e_n\|_{X^{0, b_0}_\tau} \leq C_T \left( ( T_1^{\varepsilon_0} + C_{T,s_1, \delta}\tau^{\rho}) \|e_n\|_{X^{0, b_0}_\tau} + C_{T, s_1,\delta}\tau^{- {s_1 \over 2}}  \|e_n\|_{X^{0, b_0}_{\tau}}^2
+ C_{T,s_1,\delta} \tau^{-s_1} \|e_n\|_{X^{0, b_0}_\tau}^3  \right.\\
\left. + C_{T,\delta}\tau^\rho \|e_n\|_{X^{0, b_0}_\tau}^4 + C_{T, \delta}\tau^{-s_0}\|e_n\|_{X^{0, b_0}_\tau}^5 + \tau^{s_0\over2}\right)
\end{multline*}
for some $\rho>0$. We then choose $T_{1}$ and $\tau$ small enough with respect to $C_T$, so that \\$C_{T}( T_1^{\varepsilon_0} + C_{T,s_1 ,\delta} \tau^{\rho})\leq\frac12$.
This yields
\begin{multline*}
\qquad\|e_{n}\|_{X^{0, b_0}_\tau} \leq C_{T,s_1,\delta}\tau^{-{s_1\over 2}} \|e_n\|_{X^{0, b_0}_\tau}^2 +C_{T, s_1,\delta} \tau^{-s_1} \|e_n\|_{X^{0,b_0}_\tau}^3 \\
+ C_{T,\delta} \tau^\rho \|e_n\|_{X^{0, b_0}_\tau}^4 + C_{T,\delta} \tau^{-s_0} \|e_n\|_{X^{0, b_0}_\tau}^5 + C_T \tau^{s_0\over 2}.\qquad
\end{multline*}
We next deduce that for $\tau$ sufficiently small
$$
\|e_n\|_{X^{0,b_0}_\tau}\lesssim\tau^{s_0\over 2}.
$$
Again, from the above estimates, one can easily deduce the existence of a fixed point for \eqref{m} on $[0, T_{1}]$.
This proves the desired estimate for $0\leq n\tau\leq T_1$ (which is uniform with respect to $T_1$ satisfying $ C_T(T_1^{\varepsilon_0}+C_{T,s_1,\delta}\tau^\rho)\leq\frac12$). Since the choice of $T_1$ depends only on $T$, we can then reiterate the argument on $T_1 \leq n\tau \leq 2T_1$ and so on, to get finally the estimate for $0\leq n\tau\leq T$ for $\tau$ sufficiently small (note that the needed smallness on $\tau$ depends on $T$ as usual in nonlinear problems).

\section{Proof of Theorem~\ref{mainthm}}\label{sectionmainthm}

To estimate the error $\Vert u(t_n)-u_n\Vert_{L^2}$, we just use that
$$
\Vert u(t_n)-u_n\Vert_{L^2}\leq\Vert u(t_n)-u^\tau(t_n)\Vert_{L^2}+\Vert u^\tau(t_n)-u_n\Vert_{L^2}\leq\Vert u-u^\tau\Vert_{L^\infty([0,T], L^2)}+\Vert e_n\Vert_{l^\infty_{\tau}(0 \leq n\tau \leq T, L^2)}
$$
Next, we use the embeddings \eqref{d} and \eqref{y} combined  with Theorem~\ref{theou-utau} and Proposition~\ref{propen} to get  that
$$
\Vert u(t_n)-u_n\Vert_{L^2} \leq C_{T} \tau^{s_{0} \over 2}.
$$
This concludes the proof of \eqref{c}.
\end{proof}

\section{Proof of Theorem~\ref{disckeythm}}\label{sectionproof}

Since the proof of Theorem~\ref{disckeythm} is rather long, we arrange this section as follows: we first give some technical lemmas in order to get some crucial frequency localized bilinear estimate (Lemma~\ref{lemmal4L4}). Next we will show that \eqref{z} is a consequence of these estimates. They will also be useful to obtain \eqref{r} and \eqref{rbis}. These will be done below.

This section is an adaptation to our discrete framework of the proof of the continuous case in \cite{Bourg, Burq}. We follow rather closely the steps of the proof in \cite{Burq}. The extension \eqref{rbis} will be obtained by a refinement of the proof of \eqref{r}.

Let us define for $\sigma\in  I_{\tau}= (-\frac\pi\tau, \frac\pi\tau]$,  and $m\in\mathbb{N}_0$,
$$
\mathrm{1}_{m}(\sigma)= \mathrm{1}_{ 2^m \leq \left<\sigma \right><2^{m+1} \cap I_{\tau}}, \qquad \sigma \in  I_{\tau}= (-\tfrac\pi\tau, \tfrac\pi\tau], \quad m\in\mathbb{N}_0,
$$
where $\mathrm{1}_{S}$ denotes the characteristic function of the set $S$. We then still denote  by $\mathrm{1}_{m}$ the $\frac{2\pi}\tau$ periodic extension of this function. Further, we define the operators $P_{m}$ acting on sequences of functions $(u_{n}(x)) $ by its Fourier transform
$$
\widetilde{P_{m} u_{n}}( \sigma, k) = \mathrm{1}_{m}(\sigma -|k|^2) \widetilde{u_{n}}(\sigma, k), \quad m \in \mathbb{N}_0.
$$
We also set
$$
P_{\leq M} u_{n}=  \sum_{0 \leq m \leq M} P_{m} u_{n}.
$$

In a similar way, we define localizations in spatial frequencies. We set
\begin{equation}\label{Qldef}
\widetilde{Q_{l}u_{n} }(\sigma, k)=  \mathrm{1}_{  2^l\leq \left<k\right> < 2^{l+1}} \widetilde{u_{n}}(\sigma, k)
\end{equation}
so that
$$
u_{n}= \sum_{l \geq 0} Q_{l} u_{n}.
$$
and
$$
Q_{\leq L} = \sum_{0 \leq l \leq L} Q_{l}.
$$

\begin{lemma}\label{lemdebut}
Let $a,b\ge 0$ be integers and $A=2^a$ and $B=2^b$. For a sequence of  functions $(f_n)\in l_\tau^2L^2$ such that $P_{\leq  b} Q_{\leq  a} f_{n}= f_{n}$, we have uniformly in $\tau \in (0, 1]$ that
$$
\Vert\Pi_\tau f_n\Vert_{l_\tau^4L^4}\lesssim A^\frac{\varepsilon}{2}B^{\frac{1}{2}}\Vert f_n\Vert_{l_\tau^2L^2},
$$
for any $\varepsilon>0$.
\end{lemma}

\begin{proof}
To prove this lemma we adapt the proof given in \cite{Bourg} from the continuous case. The presence of $\Pi_{\tau}$ is crucial to get a $\tau$ independent estimate.

It suffices to prove that
$$
\Vert\Pi_\tau f_n\Pi_\tau f_n\Vert_{l_\tau^2L^2}\lesssim A^\varepsilon B\Vert\Pi_\tau f_n\Vert_{l_\tau^2L^2}^2.
$$

Without loss of generality, we assume $\Vert\Pi_\tau f_n\Vert_{l_\tau^2L^2}=1.$ Thus, by Plancherel's theorem, it suffices to show that
$$
\left\Vert\sum\limits_{\substack{k_{11}+k_{21}=k_1\\k_{12}+k_{22}=k_2}}\int_{\sigma_1+\sigma_2=\sigma}\Pi_\tau\widetilde{f_n}(\sigma_1,k_{11},k_{12})\Pi_\tau\widetilde{f_n}(\sigma_2,k_{21},k_{22})d\sigma_1\right\Vert_{L^2l^2}\lesssim A^{\varepsilon}B.
$$
We have by assumption $\sigma_i-k_{i1}^2-k_{i2}^2\in E_B$, \, $i=1, 2$ with $E_B\subset\bigcup_{m \in \mathbb{Z} }[\frac{2m\pi}{\tau}-B,\frac{2m\pi}{\tau}+B]$. Since $ \sigma_{1}, \, \sigma\in I_{\tau}$ and $k_{ij}$ are restricted to $| k_{ij}|\leq \tau^{-{1 \over 2}}$ by the presence of $\Pi_{\tau}$, we actually have $E_{B}=\bigcup_{|m|\leq N}[\frac{2m\pi}{\tau}-B,\frac{2m\pi}{\tau}+B]$ with $N=\mathcal{O}(1)$, i.e. $|E_B|=\mathcal{O}(B)$. By Fubini's theorem, we have
$$
\left\Vert\sum\limits_{\substack{k_{11}+k_{21}=k_1\\k_{12}+k_{22}=k_2}}\int_{\sigma_1+\sigma_2=\sigma}|\Pi_\tau\widetilde{f_n}(\sigma_1,k_{11},k_{12})|^2|\Pi_\tau\widetilde{f_n}(\sigma_2,k_{21},k_{22})|^2d\sigma_1\right\Vert_{L^1l^1}=1,
$$
so by Cauchy--Schwarz, it will suffice to show that
$$
\sum\limits_{\substack{k_{11}+k_{21}=k_1\\k_{12}+k_{22}=k_2}}\int_{\Omega}1d\sigma_1\lesssim A^{2\varepsilon}B^2
$$
for $\Omega=\{\sigma_1\mid\sigma_1+\sigma_2=\sigma,~\sigma_1-k_{11}^2-k_{12}^2\in E_B,~\sigma_2-k_{21}^2-k_{22}^2\in E_B\}.$

Fix $k_1,k_2,\sigma$. To make the integral nonzero, we must have $\sigma-\sum\limits_{i,j=1}^2k_{ij}^2\in E_B$; then the integral is of size $\mathcal{O}(B)$. Thus it suffices to show that
$$
\sum 1\lesssim A^{2\varepsilon}B,
$$
where the sum extends over all $(k_{11},~k_{12},~k_{21},~k_{22})$ satisfying
\begin{equation}\label{7}
k_{11}+k_{21}=k_1,\quad k_{12}+k_{22}=k_2,\quad\sigma-\sum\limits_{i,j=1}^2k_{ij}^2\in E_B.
\end{equation}
Due to \eqref{7}, we have $\sigma-\frac12k_1^2-\frac12k_2^2-\frac12S\in E_B$, where $S=(k_{11}-k_{21})^2+(k_{12}-k_{22})^2$.

For $A^2>10B$, note that there are at most $\mathcal{O}(B)$ different values of $S$, thus we only need to show that there are $\mathcal{O}(A^{2\varepsilon})$ different vectors $(k_{11},~k_{12},~k_{21},~k_{22})$ if we fix $S$.

Since we have by definition that $S=\mathcal{O}(A^2)$, it suffices to show that the number of elements in the set $P=\{(c,d)\in\mathbb{Z}^2\mid c^2+d^2=S\}$ is $\mathcal{O}(S^{\varepsilon})$.

This is a consequence of  the divisor bound of $S$ in the ring of Gauss integers (see \cite[1.6.2]{Tao2}, for example):
$$
\#\{z\in\mathbb{Z}[i]\mid\exists z^\prime\in\mathbb{Z}[i], \text{ such that }zz^\prime=S\}=\mathcal{O}(S^{\varepsilon}),\quad\forall\varepsilon>0.
$$
From this, $|P|=\mathcal{O}(S^{\varepsilon})$ follows.

For $A^2\leq10B$, the point $(k_{11}-k_{21},k_{12}-k_{22})$ lies in a disk of radius $\mathcal{O}(A)\lesssim\mathcal{O}(B^{\frac12})$ in which we thus have $\mathcal{O}(B)$ integer points.

This ends the proof of Lemma~\ref{lemdebut}.
\end{proof}

From this we can then deduce the following result.

\begin{lemma}\label{lemcarre}
Let $(f_n)\in l_\tau^2L^2$ be an arbitrary sequence of functions for which $P_{l}  f_{n}= f_{n}$ and for which the support of $\widehat{f_{n}}$ is included in $a + \mathcal{C}_{N}$, where $\mathcal{C}_{N}= [-N, N)^2$ with  $a \in \mathbb{Z}^2$, $N \in \mathbb{Z}_{+}$ and $|a| \lesssim \tau^{-{1 \over 2}}.$ Then, we have the estimate (which is in particular uniform in $\tau$,  $a$ and $N$)
$$
\Vert\Pi_\tau f_n\Vert_{l_\tau^4L^4}\lesssim N^\frac{\varepsilon}{2} 2^{l \over 2}\Vert f_n\Vert_{l_\tau^2L^2},
$$
for any $\varepsilon>0$.
\end{lemma}

\begin{proof}
Let us set
\begin{equation}\label{defgalile}
g_n(x)= e^{-i \left<a,x\right>} e^{- i|a|^2n \tau} (\Pi_\tau f_n)(x + 2an\tau).
\end{equation}
We first observe that
$$
\|\Pi_\tau f_n\|_{l^4_\tau L^4}=\|g_n\|_{l^4_\tau L^4}.
$$
Moreover, we have that
$$
\widetilde{g_n}(\sigma,k)=\widetilde{\Pi_\tau f_n}\left(\sigma-|a|^2 + 2\left<a,k+a\right>,a+k\right).
$$
From this expression, we get that $\widehat{g_{n}}$ is supported in $\mathcal{C}_{N}$ and since
$$
\sigma-|a|^2+2{\left<a,k+a\right>}-|a+k|^2=\sigma- |k|^2,
$$
we get from the assumption that $P_lf_n= f_n$  that
$$
P_lg_n= g_n.
$$
Moreover, since $|a|\lesssim\tau^{-{1\over2}}$, we also have that $\Pi_{c\tau} g_n= g_n$ for some $c>0$ sufficiently small. Consequently by using Lemma~\ref{lemdebut} for $g_n$, we get that
$$
\| \Pi_{\tau} f_{n }\|_{l^4_{\tau}L^4}=\|g_{n}\|_{l^4_{\tau}L^4} = \|\Pi_{c\tau} g_{n}\|_{l^4_{\tau}L^4} \lesssim  2^{l \over 2} N^{\varepsilon \over 2} \|g_{n}\|_{l^2_{\tau} L^2}.
$$
Thanks to \eqref{defgalile}, we finally observe that
$$
\|g_{n}\|_{l^2_{\tau} L^2} =  \|f_{n}\|_{l^2_{\tau} L^2}
$$
to finally conclude the proof.
\end{proof}

Our next step will be to show the following result.
\begin{lemma}\label{lemmal4L4}
For any $\varepsilon>0$, $b>\frac{1}{2}$, $f_n\in X_\tau^{0,b}$ and $g_n\in X_\tau^{0,b}$ such that the space Fourier transforms $\widehat{f_{n}}$ and $\widehat{g_{n}}$ are supported in $ N_{1} \leq |k| \leq 2 N_{1}$ and $ N_{2} \leq |k| \leq 2 N_{2}$, respectively, with $N_{1} \leq N_{2}$, we have
$$
\Vert \Pi_\tau f_n\Pi_\tau g_n\Vert_{l_\tau^2L^2}\lesssim N_{1}^\varepsilon \Vert  f_n\Vert_{X_\tau^{0,b}}\Vert g_n\Vert_{X_\tau^{0,b}}\lesssim N_{2}^\varepsilon\Vert f_n\Vert_{X_\tau^{0,b}}\Vert g_n\Vert_{X_\tau^{0,b}}.
$$
\end{lemma}

\begin{proof}
The second inequality is trivial since $N_1\leq N_2$. For the first inequality, note that we can always assume that $\Pi_{\tau} f_{n}=f_{n}$, $\Pi_{\tau} g_{n}=g_{n}$. We then expand
$$
f_{n}= \sum_{p \geq 0} P_{p} f_{n}, \quad g_{n}= \sum_{q} P_{q} g_{n}
$$
so that
\begin{equation}\label{proofbilin1}
\Vert \Pi_\tau f_n\Pi_\tau g_n\Vert_{l_\tau^2L^2}\lesssim \sum_{p,q} \| P_{p} f_{n} P_{q} g_{n}\|_{l^2_{\tau}L^2}.
\end{equation}
Moreover, we observe that for any $b$ and for every $u_{n}$ such that $\Pi_{\tau}u_{n}=u_{n}$, we have
\begin{equation}\label{ortholit}
\|\Pi_{\tau} u_{n}\|_{X^{0, b}_{\tau}}^2\sim \sum_{p} 2^{2 p b}\| P_{p} u_{n}\|_{l^2_{\tau}L^2}^2.
\end{equation}

We shall next expand
$$
P_{q} g_{n}= \sum_{a} R_{a, N_{1}} g_{n},
$$
where $R_{a, N_{1}}$ is a projection associated to a partition of $\mathbb{Z}^2$ into cubes of  side size $N_{1}$ of the form $a + \mathcal{C}_{N_{1}}$, that is to say
$$
\widetilde{R_{a, N_{1}} g_{n}}(\sigma, k) = \mathrm{1}_{a+ \mathcal{C}_{N_{1}}}(k) \widetilde{g_{n}}(\sigma, k).
$$
Note that the sum is actually finite, with  $|a| \lesssim \tau^{-{1 \over 2}}$ by the assumption on $g_{n}$ and that, by orthogonality of the terms, we have  that
\begin{equation}\label{quasiortho}
\|P_{q} g_{n}\|_{l^2_{\tau}L^2}^2=  \sum_{a} \|R_{a, N_{1}} g_{n}\|_{l^2_{\tau}L^2}^2.
\end{equation}
We shall use the expansion
$$
\| P_{p} f_{n} P_{q} g_{n}\|_{l^2_{\tau}L^2}=\left\| \sum_{a}P_{p}f_{n}R_{a, N_{1}} g_{n} \right\|_{l^2_{\tau}L^2}.
$$
Since  the space Fourier transform of  $P_{p}f_{n}$ is supported in $B(0, N_{1})$ by assumption, the one of $P_{p}f_{n}R_{a, N_{1}} g_{n}$ is  supported in  $a+ \mathcal{C}_{10N_1}$. We thus get that the terms in the above sum are quasi-orthogonal, which yields
$$
\| P_{p} f_{n} P_{q} g_{n}\|_{l^2_{\tau}L^2} \lesssim\left(\sum_{a} \|P_{p}f_{n}R_{a, N_{1}} g_{n} \|_{l^2_{\tau}L^2}^2 \right)^{1 \over 2}.
$$
For each term inside the sum, we can then use Cauchy--Schwarz, and Lemmas~\ref{lemdebut} and \ref{lemcarre} to get
$$
\| P_{p} f_{n} P_{q} g_{n}\|_{l^2_{\tau}L^2} \lesssim\left(\sum_{a} \|P_{p}f_{n}\|_{l^4_{\tau}L^4}^2 \|R_{a, N_{1}} g_{n} \|_{l^4_{\tau}L^4}^2 \right)^{1 \over 2}
\lesssim  \left(2^{p} N_{1}^{\varepsilon} \|P_{p}f_{n}\|_{l^2_\tau L^2}^2 \sum_{a} 2^q N_{1}^\varepsilon  \|R_{a, N_{1}} g_{n} \|_{l^2_{\tau}L^2}^2 \right)^{1 \over 2}.
$$
Therefore, from \eqref{quasiortho}, we get
$$
\| P_{p} f_{n} P_{q} g_{n}\|_{l^2_{\tau}L^2} \lesssim  2^{p\over 2} N_{1}^{\varepsilon} \|P_{p}f_{n}\|_{l^2_\tau L^2}  2^{q \over 2}\|P_{q}g_{n} \|_{l^2_{\tau}L^2} .
$$
To conclude, we use \eqref{proofbilin1} and Cauchy--Schwarz
$$
\Vert \Pi_\tau f_n\Pi_\tau g_n\Vert_{l_\tau^2L^2}\lesssim N_{1}^{\varepsilon} \sum_{p\geq 0,q \geq0}  2^{p ({1 \over 2} -b)}  2^{pb} \|P_{p}f_{n}\|_{l^2_\tau L^2}  2^{q({1 \over 2}- b)}2^{qb}\|P_{q}g_{n} \|_{l^2_{\tau}L^2}\lesssim N_{1}^{\varepsilon} \| f_{n}\|_{X^{0, b}_{\tau}}  \| g_{n}\|_{X^{0, b}_{\tau}}
$$
for $b>\frac12.$
Note that we have also used \eqref{ortholit} to get the last estimate.
\end{proof}

We can now deduce \eqref{z}.

\begin{corollary}\label{corStrich}
For every $\varepsilon>0$ and every $b>\frac12$, we have
$$
\| \Pi_\tau f_n\|_{l^4_{\tau}L^4} \lesssim \| f_n\|_{X^{\varepsilon, b}_{\tau}}.
$$
\end{corollary}

\begin{proof}
We use the Littlewood--Paley decomposition in space. By using \eqref{Qldef}, we write
$$
\|\Pi_\tau f_n\|_{l^4_{\tau}L^4} \leq \sum_{l\geq 0} \| Q_{l} \Pi_\tau f_n\|_{l^4_{\tau}L^4}.
$$
From Lemma~\ref{lemmal4L4}, we thus get
$$
\| \Pi_{\tau} f_{n}\|_{l^4_{\tau}L^4} \lesssim \sum_{l \geq 0} 2^{l \varepsilon \over 2} \|Q_{l}f_{n}\|_{X_{\tau}^{0, b}}.
$$
From Cauchy--Schwarz and a standard characterization of the $H^s$ norm, we get
$$
\| \Pi_{\tau} f_{n}\|_{l^4_{\tau}L^4} \lesssim \left( \sum_{l \geq 0}  2^{ 2l \varepsilon } \|Q_{l}f_{n}\|_{X_{\tau}^{0, b}}^2\right)^{1 \over 2}\lesssim \|f_{n}\|_{X_\tau^{\varepsilon, b}}.
$$
This concludes the proof.
\end{proof}

We now turn to the proof of \eqref{r}.
\subsection*{ Proof of  \eqref{r}}

Take $u_{0,n}$ to be any function in $X_\tau^{-s,b}$. By duality, the problem is equivalent to
$$
I=\left|\tau\sum\limits_{n}\int_{\mathbb{T}^2} \overline{\Pi_{\tau}u_{0,n}} \Pi_\tau u_{1,n} \overline{\Pi_\tau u_{2,n}}\Pi_\tau u_{3,n}\,dx\right|\lesssim\Vert u_{0,n}\Vert_{X_\tau^{-s,b}}\Vert u_{1,n}\Vert_{X_\tau^{s,b}}\Vert u_{2,n}\Vert_{X_\tau^{s,b}}\Vert u_{3,n}\Vert_{X_\tau^{s,b}}.
$$
We shall use again Littlewood--Paley decompositions. We set $N=(N_0,N_1,N_2,N_3),~L=(L_0,L_1,L_2,L_3)$, where $N_{j}$ and $L_{j}$ are dyadic of the form
$N_j=2^{n_j},L_j=2^{l_j},~j=0,1,2,3$. We split
$$
\Pi_{\tau}u_{j,n}=\sum\limits_{l_j,n_j\in\mathbb{N}}u_{j, n}^{L_jN_j}, \qquad u_{j,n}^{L_jN_j} = P_{l_{j}} Q_{n_{j}} u_{j,n}, \quad j=0,1,2,3.
$$
We shall use that for any $\sigma$, $\beta$, we have
\begin{equation}\label{equivlittle}
\| \Pi_{\tau}u_{j,n} \|_{X^{\sigma, \beta}_{\tau}} \sim\left( \sum\limits_{l_j,n_j\in\mathbb{N}} N_{j}^{2 \sigma} L_{j}^{2 \beta} \|u_{j, n}^{L_jN_j}\|_{l^2_{\tau}L^2}^2
\right)^{1 \over 2}.
\end{equation}
For notational convenience, we thus set
$$
c_0(L_0, N_0)^2\sim N_j^{-2s} L_j^{2b} \|u_{0,n}^{L_0N_0}\|_{l^2_{\tau}L^2}^2, \qquad
c_j(L_j, N_j)^2 \sim N_j^{2s} L_j^{2b} \|u_{j, n}^{L_jN_j}\|_{l^2_{\tau}L^2}^2, \quad j=1,\,2,\,3
$$
so that
$$
\sum_{L_0, N_0} c_0(L_0, N_0)^2\sim \Vert u_{0,n}\Vert^2_{X_\tau^{-s,b}}, \qquad
\sum_{L_j, N_j} c_j(L_j, N_j)^2\sim \Vert u_{j,n}\Vert^2_{X_\tau^{s,b}}, \quad j=1,\,2,\,3.
$$

Let us also set
\begin{equation}\label{defILN}
I(L,N)=\left|\tau\sum_n\int_{\mathbb{T}^2} \overline{u_{0,n}^{L_0N_0}}u_{1,n}^{L_1N_1}
\overline{u_{2,n}^{L_2N_2}}u_{3,n}^{L_3N_3} dx\right|
\end{equation}
so that
$$
I \leq \sum_{L, N} I(L,N).
$$
Note that, by properties of the support in Fourier space of the involved functions, $I(L, N)$ vanishes unless $N_0\lesssim N_1+N_2+N_3.$ To estimate the above sum, we can restrict ourselves without loss of generality to the case that  $N_1\geq N_2\geq N_3$. This yields in particular that  $N_1+N_2+N_3\lesssim N_1$, and therefore, we must have $ N_{0} \lesssim N_{1}$ so that $I(L, N)$ does not vanish.

By Cauchy--Schwarz and Lemma~\ref{lemmal4L4}, for any $\varepsilon_1>0$ to be chosen small enough, we have
\begin{equation}\label{2}
\begin{aligned}
I(L,N)&\lesssim\Vert u_{0,n}^{L_0N_0}u_{2,n}^{L_2N_2}\Vert_{l_\tau^2L^2}\Vert u_{1,n}^{L_1N_1}u_{3,n}^{L_3N_3}\Vert_{l_\tau^2L^2}\\
&\lesssim (N_2N_3)^{\varepsilon_1}\Vert u_{0,n}^{L_0N_0}\Vert_{X_\tau^{0,\frac{1}{2}+\varepsilon_1}}\Vert u_{1,n}^{L_1N_1}\Vert_{X_\tau^{0,\frac{1}{2}+\varepsilon_1}}\Vert u_{2,n}^{L_2N_2}\Vert_{X_\tau^{0,\frac{1}{2}+\varepsilon_1}}\Vert u_{3,n}^{L_3N_3}\Vert_{X_\tau^{0,\frac{1}{2}+\varepsilon_1}}\\
&\lesssim(L_0L_1L_2L_3)^{\frac{1}{2}+\varepsilon_1-b}\dfrac{N_0^s}{(N_1N_2N_3)^s}(N_2N_3)^{\varepsilon_1}\prod\limits_{j=0}^3c_j(L_j,N_j).
\end{aligned}
\end{equation}
We shall then perform another estimate of $I(L,N)$ which behaves in a better way with respect to the $L$ localizations. We first observe that, by interpolation between the continuous  embeddings $X^{0, b}_{\tau} \subset l^\infty_\tau L^2$ for every $b>\frac12$ and the trivial one $X^{0,0}_\tau\subset l^2_{\tau}L^2$, we get that $X^{0,b}_\tau\subset l^4_{\tau}L^2$ for every $b>\frac14$, that is to say, there exists $C>0$ independent of $\tau \in(0, 1]$ such that
$$
\|u_{n}\|_{l^4_{\tau}L^2} \leq C   \|u_{n}\|_{X^{0, b}_{\tau}}, \quad \text{for all } (u_{n}) \in X^{0, b}_{\tau}.
$$
By using H\"{o}lder's inequality, the previous estimate with $b= {1 \over 4}+ \varepsilon_{1}$ and the Sobolev embedding $H^{1+\varepsilon_1}(\mathbb{T}^2)\subset L^{\infty}(\mathbb{T}^2)$, we have
\begin{equation}\label{3}
\begin{aligned}
I(L,N)&\lesssim\Vert u_{0,n}^{L_0N_0}\Vert_{l_\tau^4L^2}\Vert u_{1,n}^{L_1N_1}\Vert_{l_\tau^4L^2}\Vert u_{2,n}^{L_2N_2}\Vert_{l_\tau^4 L^\infty}\Vert u_{3,n}^{L_3N_3}\Vert_{l_\tau^4 L^\infty}\\
&\lesssim\Vert u_{0,n}^{L_0N_0}\Vert_{X_{\tau}^{0, {1 \over 4}+ \varepsilon_{1}}}\Vert u_{1,n}^{L_1N_1}\Vert_{X_{\tau}^{0, {1\over 4}+ \varepsilon_{1}}}\Vert u_{2,n}^{L_2N_2}\Vert_{X_{\tau}^{1+ \varepsilon_{1}, {1 \over 4} + \varepsilon_{1} }}\Vert u_{3,n}^{L_3N_3}\Vert_{X_{\tau}^{ 1 +\varepsilon_{1}, {1 \over 4}+ \varepsilon_{1}}}\\
&\lesssim(L_0L_1L_2L_3)^{\frac{1}{4}+\varepsilon_1-b}\dfrac{N_0^s}{(N_1N_2N_3)^s}
(N_2N_3)^{1+\varepsilon_{1}}\prod\limits_{j=0}^3c_j(L_j,N_j).
\end{aligned}
\end{equation}
We shall then interpolate between  \eqref{3} with strength $\theta \in (0, 1)$ and \eqref{2} with strength $(1-\theta)$. This yields
$$
I(L, N) \lesssim  \left({N_{0} \over N_{1}} \right)^s { 1 \over (L_{0}L_{1}L_{2}L_{3})^{\mu_{1}}} { 1 \over (N_{2}N_{3})^{\mu_{2}}} \prod\limits_{j=0}^3c_j(L_j,N_j)
$$
with
$$ \mu_{1}= b-  \theta(\tfrac14 + \varepsilon_{1}) - (1- \theta)(\tfrac12 + \varepsilon_{1}),
\quad \mu_{2}= s - \theta(1 +  \varepsilon_{1})  - (1- \theta) \varepsilon_{1}.
$$
We then observe that we can choose $\theta \in (0, 1)$ such that $ s>\theta > 2 - 4b$ and then  $\varepsilon_{1}$ sufficiently small to get $\mu_{1}>0$ and $\mu_{2}>0$.
With this choice, we obtain easily by Cauchy--Schwarz that
\begin{align*}
\sum\limits_{L,N}I(L,N)
&\lesssim\sum\limits_{N_0,N_1}\left(\dfrac{N_0}{N_1}\right)^s\left(\sum\limits_{L,N_2,N_3}\dfrac{1}{(N_2N_3)^{2\mu_2}(L_0L_1L_2L_3)^{2\mu_1}}\right)^{\frac{1}{2}}\\
&\qquad\left(\sum\limits_{L_0}c_0^2(L_0,N_0)\sum\limits_{L_1}c_1^2(L_1,N_1)\sum\limits_{L_2,N_2}c_2^2(L_2,N_2)\sum\limits_{L_3,N_3}c_3^2(L_3,N_3)\right)^{\frac{1}{2}}\\
&\lesssim\Vert u_{2,n}\Vert_{X_\tau^{s,b}}\Vert u_{3,n}\Vert_{X_\tau^{s,b}}\sum\limits_{N_0,N_1}\left(\dfrac{N_0}{N_1}\right)^s\left(\sum\limits_{L_0}c_0^2(L_0,N_0)\right)^{\frac{1}{2}}\left(\sum\limits_{L_1}c_1^2(L_1,N_1)\right)^{\frac{1}{2}}.
\end{align*}
To conclude, we can use again that the sum is restricted to  $N_0\lesssim N_1$, we can thus write  $N_1=2^\gamma N_0$ for  $\gamma\geq\gamma_0$, where $\gamma_0$ is a fixed integer. This yields
\begin{align*}
\sum\limits_{N_0,N_1}\left(\dfrac{N_0}{N_1}\right)^s&\left(\sum\limits_{L_0}c_0^2(L_0,N_0)\right)^{\frac{1}{2}}\left(\sum\limits_{L_1}c_1^2(L_1,N_1)\right)^{\frac{1}{2}}\\
&=\sum\limits_{\gamma\geq\gamma_0}2^{-s\gamma}\sum\limits_{N_0}\left(\sum\limits_{L_0}c_0^2(L_0,N_0)\right)^{\frac{1}{2}}\left(\sum\limits_{L_1}c_1^2(L_1,2^\gamma N_0)\right)^{\frac{1}{2}}\\
&\lesssim\sum\limits_{\gamma\geq\gamma_0}2^{-s\gamma}\left(\sum\limits_{L_0,N_0}c_0^2(L_0,N_0)\right)^{\frac{1}{2}}\left(\sum\limits_{L_1,N_0}c_1^2(L_1,2^\gamma N_0)\right)^{\frac{1}{2}}\\
&\lesssim\Vert u_{0,n}\Vert_{X_\tau^{-s,b}}\Vert u_{1,n}\Vert_{X_\tau^{s,b}},
\end{align*}
where we have used Cauchy--Schwarz again to pass from the second line to the third line.
This ends the proof of \eqref{r}.

\subsection*{Proof of \eqref{rbis}}
We follow the same strategy as above and use the same notations. For any  $(u_{0,n}) \in X_\tau^{0,b}$, it is equivalent to prove that
$$
I=\left|\tau\sum\limits_{n}\int_{\mathbb{T}^2}\overline{ \Pi_{\tau}u_{0,n}} \Pi_\tau u_{1,n}\overline{\Pi_\tau u_{2,n}}\Pi_\tau u_{3,n}\,dx\right|\lesssim\Vert u_{0,n}\Vert_{X_\tau^{0,b}}\Vert u_{\sigma(1),n}\Vert_{X_\tau^{s,b}}\Vert u_{\sigma(2),n}\Vert_{X_\tau^{s,b}}\Vert u_{\sigma(3),n}\Vert_{X_\tau^{0,b}}.
$$
By using again the definition \eqref{defILN}, the proof follows if we prove that
$$
S:= \sum_{\substack{L, N \\ N_{3} \leq N_{2} \leq N_{1}} }I(L, N) \lesssim
\Vert u_{0,n}\Vert_{X_\tau^{0,b}}\Vert u_{\sigma(1),n}\Vert_{X_\tau^{s,b}}\Vert u_{\sigma(2),n}\Vert_{X_\tau^{s,b}}\Vert u_{\sigma(3),n}\Vert_{X_\tau^{0,b}},
$$
where $\sigma \in \mathcal{S}_{3}$ is any permutation of $\{ 1, 2, 3 \}$. By a straightforward adaptation of the previous estimate for $I(L,N)$, we get that for $N_{3} \leq N_{2} \leq N_{1}$,
\begin{equation}\label{ILNbis}
I(L,N) \lesssim { 1 \over (L_0L_1L_2L_3)^{\mu_1}}{ 1\over (N_2N_3)^{\mu_2}}\| u_{0,n}^{L_0, N_0} \|_{X^{0,b}_\tau}\| u_{2,n}^{L_2, N_2} \|_{X^{s, b}_\tau}\| u_{3,n}^{L_3,N_3}\|_{X^{s,b}_\tau}\| u_{1,n}^{L_1, N_1} \|_{X^{0,b}_\tau}
\end{equation}
with
$$
\mu_1= b-\theta(\tfrac14 + \varepsilon_1)-(1- \theta)(\tfrac12 + \varepsilon_1),\quad \mu_2= s-\theta(1+ \varepsilon_1) - (1-\theta)\varepsilon_1
$$
and hence we can choose $\theta $ and $\varepsilon_1$ as before to get $\mu_1>0$, $\mu_2>0$.

Let us recall that $I(L,N)$ vanishes unless $N_0\lesssim N_1$. We shall distinguish two subcases to estimate $S$.

Let $N_2\leq N_1/4$. Then, since $I(L,N)$ vanishes unless the sum of the spatial frequencies of the involved function is $0$, which is to say, $-k_0+k_1-k_2+k_3=0$, we also get that $N_0\geq N_1/2$. We can thus write $N_0= 2^\gamma N_1$, where $\gamma_0\leq \gamma \leq \gamma_1$ with two fixed values $\gamma_0$ and $\gamma_1$.
Let us set
$$
S_1:= \sum_{\substack{L,N \\ N_3\leq N_2\leq N_1/4}}I(L, N).
$$
We thus get that
$$
S_1\lesssim \sum_{\gamma_0 \leq \gamma \leq \gamma_1}\sum_{L, N_1, N_2, N_3}{ 1 \over (L_0L_1L_2L_3)^{\mu_1}}{ 1 \over (N_2N_3)^{\mu_2}} \| u_{0,n}^{L_0, 2^\gamma N_1} \|_{X^{0, b}_\tau}\| u_{2,n}^{L_2, N_2} \|_{X^{s, b}_\tau}\| u_{3,n}^{L_3, N_3} \|_{X^{s, b}_\tau}\| u_{1,n}^{L_1, N_1} \|_{X^{0, b}_\tau}
$$
and hence by Cauchy--Schwarz that
$$
S_1\lesssim \|u_{0,n}\|_{X^{0,b}_\tau} \|u_{2,n}\|_{X^{s,b}_\tau} \|u_{3,n}\|_{X^{s,b}_\tau}\|u_{1,n}\|_{X^{0,b}_\tau}.
$$
Note that this is the sharpest estimate since we do not lose space derivatives on the function which has the highest frequencies in the sum.
We can also use that $N_2 \lesssim  N_1$ to get that
\begin{multline*}
\| u_{2,n}^{L_2, N_2} \|_{X^{s,b}_\tau} \| u_{1,n}^{L_1, N_1} \|_{X^{0, b}_\tau}
\lesssim N_2^s \| u_{2,n}^{L_2, N_2} \|_{X^{0,b}_\tau}\| u_{1,n}^{L_1, N_1} \|_{X^{0,b}_\tau}
\lesssim N_1^s \| u_{2,n}^{L_2, N_2} \|_{X^{0,b}_\tau}\| u_{1,n}^{L_1, N_1} \|_{X^{0,b}_\tau}  \\
\lesssim \| u_{2,n}^{L_2, N_2} \|_{X^{0, b}_\tau}\| u_{1,n}^{L_1, N_1} \|_{X^{s,b}_\tau}
\end{multline*}
and hence deduce from \eqref{ILNbis} that $S_1$ can be also estimated by
$$
S_1 \lesssim \|u_{0,n}\|_{X^{0, b}_\tau} \|u_{2,n}\|_{X^{0,b}_\tau} \|u_{3,n}\|_{X^{s,b}_\tau}\|u_{1,n}\|_{X^{s,b}_\tau}.
$$
By using that $N_3\leq N_1$ and a similar argument, we thus actually obtain that
$$
S_1\leq \Vert u_{0,n}\Vert_{X_\tau^{0,b}}\Vert u_{\sigma(1),n}\Vert_{X_\tau^{s,b}}\Vert u_{\sigma(2),n}\Vert_{X_\tau^{s,b}}\Vert u_{\sigma(3),n}\Vert_{X_\tau^{0,b}}
$$
for any $\sigma \in \mathcal{S}_{3}.$

It remains to study
$$
S_2:= \sum_{\substack{L, N \\ N_3 \leq N_2 \leq N_1 \\ N_2 \geq N_1/4}}I(L, N).
$$
By using \eqref{ILNbis}, we now deduce from  $N_2\geq N_1/4$ and since we still have $N_0\lesssim N_1$ that
$$
I(L, N) \lesssim   { 1 \over (L_0L_1L_2L_3)^{\mu_1}} { 1 \over (N_0N_1N_2N_3)^{\mu_2 \over 3}}  \| u_{0,n}^{L_0, N_0} \|_{X^{0, b}_\tau}\| u_{2,n}^{L_{2}, N_{2}} \|_{X^{s, b}_\tau}\| u_{3,n}^{L_3, N_3} \|_{X^{s, b}_\tau}\| u_{1,n}^{L_1, N_1} \|_{X^{0, b}_\tau}.
$$
Consequently, we directly deduce from Cauchy--Schwarz that
$$
S_2 \lesssim \|u_{0,n}\|_{X^{0, b}_\tau} \|u_{2, n}\|_{X^{s,b}_\tau} \|u_{3,n}\|_{X^{s,b}_\tau} \|u_{1,n}\|_{X^{0,b}_\tau}.
$$
From the same observation as above, since $N_{2} \leq N_{1}$ and $N_3\leq N_1$, we then obtain that
$$
S_2 \lesssim \|u_{0,n}\|_{X^{0,b}_\tau}\Vert u_{\sigma(1),n}\Vert_{X_\tau^{s,b}}\Vert u_{\sigma(2),n}\Vert_{X_\tau^{s,b}}\Vert u_{\sigma(3),n}\Vert_{X_\tau^{0,b}}
$$
for any $\sigma \in \mathcal{S}_3.$

Since $S= S_1+ S_2$, this ends the proof of \eqref{rbis}.

\section{Numerical experiments}

\begin{figure}[b]
\begin{center}
\subfigure[]{\includegraphics[width=0.48\textwidth]{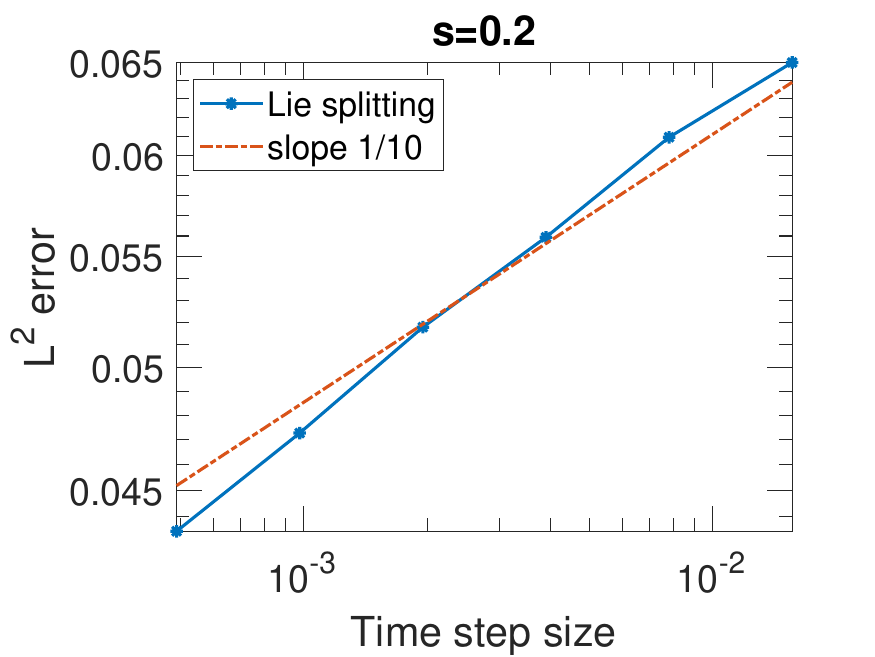}}
\subfigure[]{\includegraphics[width=0.48\textwidth]{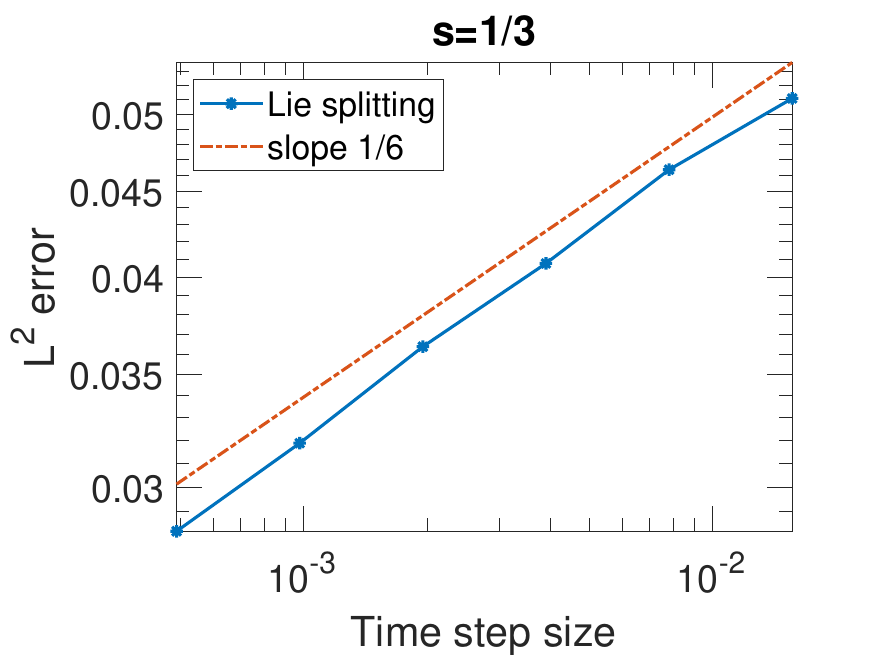}}
\subfigure[]{\includegraphics[width=0.48\textwidth]{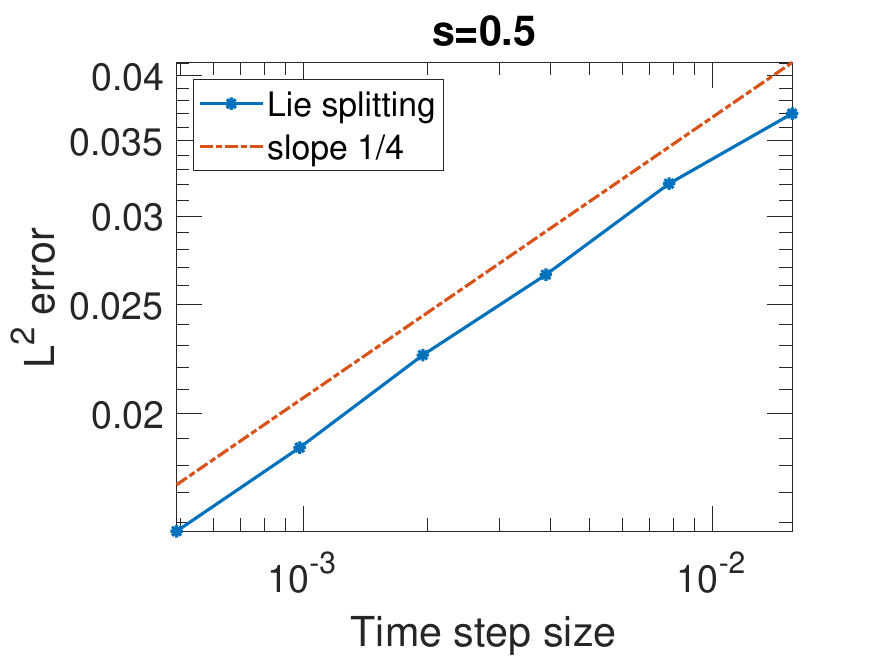}}
\subfigure[]{\includegraphics[width=0.48\textwidth]{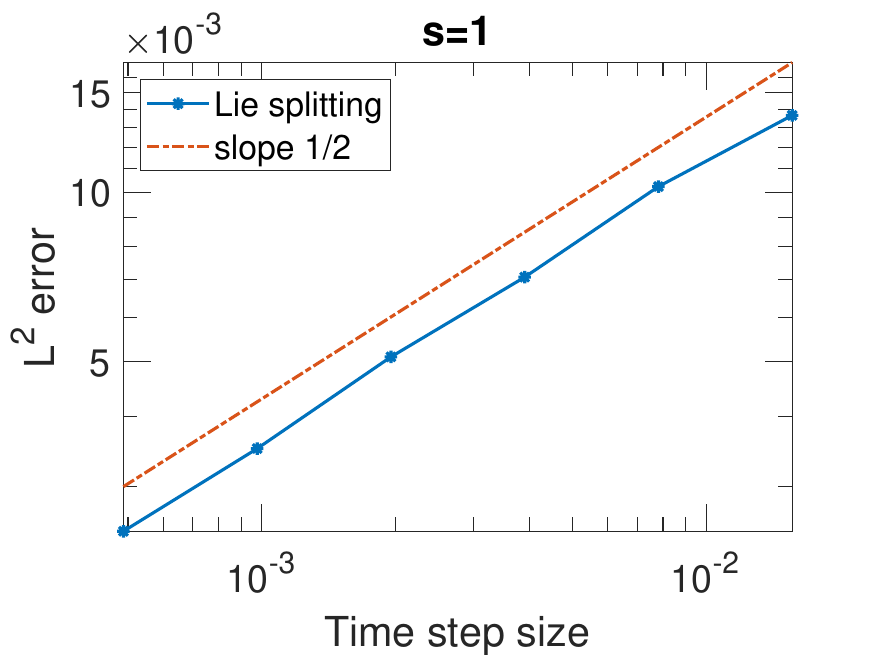}}
\end{center}
\caption{$L^2$ error of the filtered  Lie splitting scheme for rough initial data $u_0\in H^s$. \qquad
(a) $s=0.2$; \quad (b) $s={1/3}$; \quad (c) $s={0.5}$; \quad (d) $s=1$.\label{fig:figure1}}
\end{figure}

In this section, we numerically illustrate our main result (Theorem~\ref{mainthm}). We display the convergence order of the filtered Lie splitting method with rough initial data. In Figure 1 we consider the periodic NLS \eqref{n}, discretized with the filtered Lie splitting method \eqref{0} and initial data
$$
u_0=\sum\limits_{k\in\mathbb{Z}^2}\langle k\rangle^{-(s+1)}\tilde{g}_ke^{i\langle k,x\rangle}\in H^s,
$$
with $s=0.2,~1/3,~0.5$ and 1, where $\tilde{g}_k$ are random variables which are uniformly distributed in the square $[-1,1]+i[-1,1]$. We employ a standard Fourier pseudospectral method for the discretization in space and we choose as largest Fourier mode $K=(2^{12},2^{12})$, i.e., the spatial mesh size $\Delta x=0.0015$ for both of the two space directions. We normalize the $L^2$ norm of the initial data by 0.1. As a reference solution, we use the filtered Lie splitting method with $K$ spatial points and a very small time step size $\tau=2^{-22}$. We choose $T=1$ to be the final time.

From Figure~\ref{fig:figure1}, we can clearly conclude that our numerical experiments confirm the convergence rate of order $\mathcal{O}(\tau^{\frac{s}{2}})$ for solutions in $H^s$ (see Theorem~\ref{mainthm}) with $s=0.2,~1/3,~0.5$ and 1.

\begin{figure}[t]
\begin{center}
\subfigure[]{\includegraphics[width=0.48\textwidth]{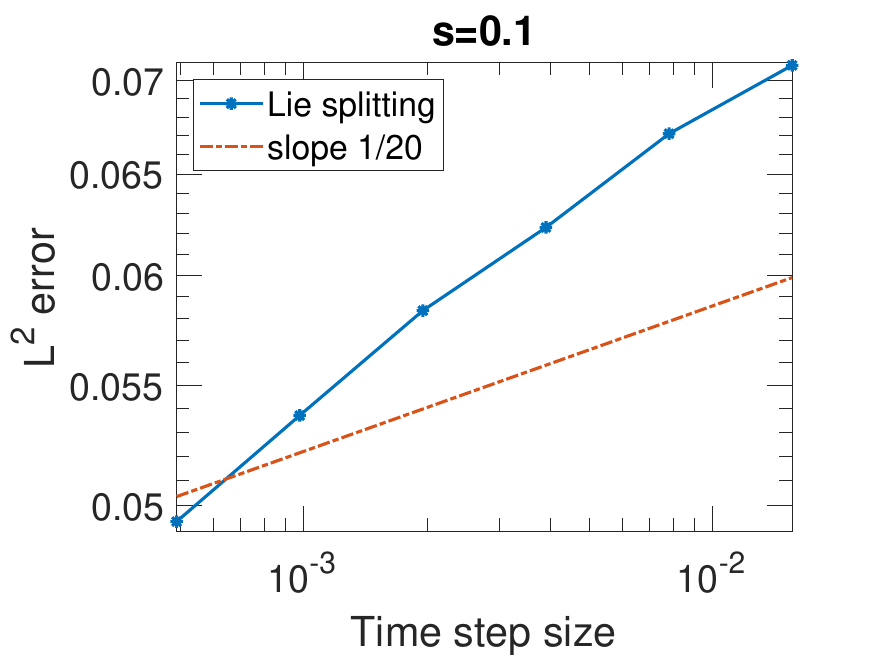}}
\subfigure[]{\includegraphics[width=0.48\textwidth]{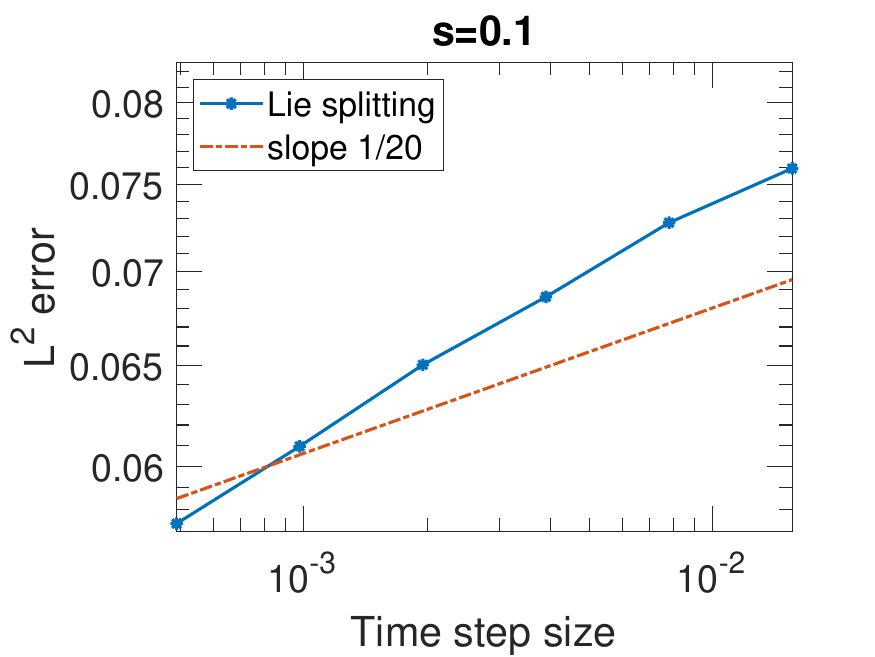}}
\end{center}
\caption{$L^2$ error of the filtered Lie splitting scheme for rough initial data $u_0\in H^{0.1}$ with different reference solutions. (a) Reference solution with largest Fourier mode $K=(2^{10},2^{10})$, spatial mesh size $\Delta x=0.0061$, and time step size $\tau=2^{-18}$; (b) Reference solution with largest Fourier mode $K=(2^{12}, 2^{12})$, spatial mesh size $\Delta x=0.0015$, and time step size $\tau=2^{-22}$.\label{fig:figure2}}
\end{figure}

We also did some experiments for very small $s$, as for example, $s=0.1$. In this case, the convergence is very slow. To obtain an accurate reference solution for small $s$, one would need an ever increasing number of Fourier modes, which is beyond the capabilities of our computers. In Figure~\ref{fig:figure2}, we demonstrate that the proven order of convergence only shows up for sufficiently high spatial resolution.

{}

\begin{thebibliography}{00}
\bibitem{BLW} G. Bai, B. Li, Y. Wu, \textit{A constructive low-regularity integrator for the one-dimensional cubic non-linear Schr\"{o}\-dinger equation under the Neumann boundary condition}, IMA J. Numer. Anal. 43: 3243--3281 (2022).
\bibitem{Bourg} J. Bourgain, \textit{Fourier transform restriction phenomena for certain lattice subsets and applications to nonlinear evolution equations. Part I: Schr\"{o}dinger equations}, Geom. Funct. Anal. 3: 209--262 (1993).
\bibitem{Burq} N. Burq, P. G\'{e}rard, N. Tzvetkov, \textit{Bilinear eigenfunction estimates and the nonlinear Schr\"{o}dinger equation on surfaces}, Invent. Math., 159: 187--223 (2005).
\bibitem{Pavlovic} D. De Silva, N.  Pavlovic, G.  Staffilani, N. Tzirakis, \textit{Global well-posedness for a periodic nonlinear Schr\"{o}dinger equation in 1D and 2D}, Discrete Contin. Dyn. Syst. 19(1): 37--65 (2007).
\bibitem{ESS16} J. Eilinghoff, R. Schnaubelt, K. Schratz, \textit{Fractional error estimates of splitting schemes for the nonlinear Schr\"{o}dinger equation}, J. Math. Anal. Appl. 442: 740--760 (2016).
\bibitem{Faou12} E. Faou, \textit{Geometric Numerical Integration and Schr\"{o}dinger Equations}, European Math. Soc. Publishing House, Z\"{u}rich 2012.
\bibitem{Ignat11} L. I. Ignat, \textit{A splitting method for the nonlinear Schr\"{o}dinger equation}, J. Differential Equations 250: 3022--3046 (2011).
\bibitem{Lubich08} C. Lubich, \textit{On splitting methods for Schr\"{o}dinger-Poisson and cubic nonlinear Schr\"{o}dinger equations}, Math. Comp. 77: 2141--2153 (2008).
\bibitem{Mus} C. Muscalu, W. Schlag, \textit{Classical and multilinear harmonic analysis}, Cambridge University Press, Cambridge, 2013.
\bibitem{ORS1} A. Ostermann, F. Rousset, K. Schratz, \textit{Error estimates of a Fourier integrator for the cubic Schr\"odinger equation at low regularity}. Found. Comput. Math. 21(3): 725--765 (2021).
\bibitem{Ost} A. Ostermann, F. Rousset, K. Schratz, \textit{Fourier integrator for periodic NLS: low regularity estimates via Bourgain spaces}, J. Eur. Math. Soc. 25(10): 3913--3952 (2023).
\bibitem{Ost1} A. Ostermann, F. Rousset, K. Schratz, \textit{Error estimates at low regularity of splitting schemes for NLS}, Math. Comp. 91: 169--182 (2022).
\bibitem{Rou} F. Rousset, K. Schratz, \textit{Convergence error estimates at low regularity for time discretizations of KdV}, Pure and Applied Analysis, 4(1): 127--152 (2022).
\bibitem{Tao} T. Tao, \textit{Nonlinear dispersive equations: local and global analysis}, Amer. Math. Soc., Providence RI, 2006.
\bibitem{Tao2} T. Tao, \textit{Poincar\'{e}'s legacies, part I: pages from year two of a mathematical blog}, Amer. Math. Soc., 2009.
\bibitem{Wu} Y. Wu, \textit{A modified splitting method for the cubic nonlinear Schr\"{o}dinger equation}, arXiv: 2212.09301.
\bibitem{WY22} Y. Wu, F. Yao, \textit{A first-order Fourier integrator for the nonlinear Schr\"{o}dinger equation on $\mathbb{T}$ without loss of regularity}, Math. Comp. 91: 1213--1235 (2022).
\end{thebibliography}
\end{document}